\documentclass[12pt]{article}
\usepackage{url}
\usepackage{epic, eepic, subfigure, floatflt}
\usepackage{amsfonts, epsfig, latexsym, amsmath}
\begin{document}
\newcommand{\R}{\ensuremath{\mathbb{R}}}
\newcommand{\N}{\ensuremath{\mathbb{N}}}
\newcommand{\Q}{\ensuremath{\mathbb{Q}}}
\newcommand{\C}{\ensuremath{\mathbb{C}}}
\newcommand{\Z}{\ensuremath{\mathbb{Z}}}
\newcommand{\T}{\ensuremath{\mathbb{T}}}
\newtheorem{theorem}{Theorem}[section]
\newtheorem{definition}[theorem]{Definition}
\newtheorem{conjecture}[theorem]{Conjecture}
\newtheorem{corollary}[theorem]{Corollary}
\newtheorem{lemma}[theorem]{Lemma}
\newtheorem{claim}[theorem]{Claim}
\newtheorem{remark}[theorem]{Remark}
\newtheorem{proposition}[theorem]{Proposition}
\newcommand{\qed}{\hfill $\Box$ }
\newcommand{\proof}{\noindent{\bf Proof.}\ \ }
\newcommand{\sketchproof}{\noindent{\bf Sketch of proof.}\ \ }

\title{\Large{\bf Zigzag Structures of Simple Two-faced Polyhedra}}

\author{Michel DEZA\footnote{Michel.Deza@ens.fr}\\
        \normalsize   LIGA, ENS/CNRS, Paris and\\
        \normalsize  Institute of Statistical Mathematics, Tokyo\\
\and
        Mathieu DUTOUR\footnote{Mathieu.Dutour@ens.fr}\\
        \normalsize  LIGA, ENS/CNRS, Paris and\\
        \normalsize  The Hebrew University, Jerusalem
\footnote{Research financed by EC's IHRP Programme, within the Research Training Network ``Algebraic Combinatorics in Europe,'' grant HPRN-CT-2001-00272.}\\
}

%       \
%\and
%       Mikhail SHTOGRIN\\
%       \normalsize   Steklov Mathematical Institute,\\
%       \normalsize   Moscow\\

\date{\small \today}
 
\maketitle

\baselineskip=20pt
\parindent=1cm

\begin{abstract}
\noindent

A {\em zigzag} in a plane graph is a circuit of edges, such that any two, 
but no three, consecutive edges belong to the same face.
A {\em railroad} in a plane graph is a circuit of hexagonal faces,
such that any hexagon is adjacent to its neighbors on opposite
edges. A graph without a railroad is called {\em tight}.
We consider the zigzag and railroad structures of general
$3$-valent plane graph and, especially, of simple two-faced
polyhedra, i.e., $3$-valent $3$-polytopes with only 
$a$-gonal and $b$-gonal faces, where $3 \le a < b \le 6$; 
the main cases are $(a,b)=(3,6)$, $(4,6)$ and $(5,6)$ (the {\em fullerenes}).

We completely describe the zigzag structure for the case $(a,b)$=$(3,6)$. 
For the case $(a,b)$=$(4,6)$ we describe symmetry groups, classify
all tight graphs with simple zigzags and give the upper bound $9$
for the number of zigzags in general tight graphs.
For the remaining case $(a,b)$=$(5,6)$ we give a construction realizing a
prescribed zigzag structure.

\end{abstract}

{\small {\em Mathematics Subject Classification}. Primary 52B05, 52B10;
Secondary 05C30,\\
\indent 05C10.

}

\bigskip
\bigskip

\section{Introduction: main notions}

A {\em graph} $G$ consists of a vertex-set $V(G)$ and an 
edge-set $E(G)$, such that either one or two vertices are
assigned to each edge as its ends. A graph is said to be {\em simple}
if no edge has only one end-vertex and no two edges have identical 
end-vertices. A {\em plane} graph is a particular drawing of a 
graph in the Euclidean plane using smooth curves that cross
each other only at the vertices of the graph. A graph, which
has at least one such drawing, is called {\em planar}.
We only deal with simple connected plane graphs, 
whose vertex-sets and edge-sets are finite. 

%A {\em walk} of length $\ell$ in a graph $G$ is a 
%sequence $v_0e_1v_1 \cdots e_{\ell}v_{\ell}$, where $v_0, v_1,
%\ldots, v_{\ell}$ are vertices of $G$, $e_1, e_2, \ldots,
%e_{\ell}$ are edges of $G$, and $v_{i-1}$ and $v_i$
%are the ends of $e_i$ for $1 \le i \le \ell$. A {\em circuit} 
%is a walk without repeated edges, such that
%$v_0=v_{\ell}$. A {\em cycle} is a circuit without
%repeated vertices. 

By a {\em polyhedron} we mean a convex $3$-dimensional polytope.
The vertices and the edges of a polyhedron form 
a simple planar $3$-connected graph called {\em skeleton}. Such a graph has 
a unique drawing (up to diffeomorphisms of the Euclidean plane),
but there exist many distinct polyhedra having this graph as 
their skeleton; all such polyhedra are of the same 
{\em combinatorial type}. The group of isometries of a polyhedron 
(so-called {\em point group}) is a subgroup of the algebraic symmetry 
group of the graph.
By a theorem of Mani (\cite{M}), there is, for each combinatorial type,
at least one polyhedron, for which this inclusion becomes equality.
So, we can identify the polyhedron and its graph, as well as the algebraic
symmetry group and the point group.

% A polyhedron $P$ is called {\em semiregular}
%if the faces of $P$ are regular and the group of 
%symmetries of $P$ acts transitively on the vertices of $P$.
%The semiregular polyhedra include the 5 well-known Platonic
%solids, 13 Archimedean solids,
%and two infinite families of prisms and antiprisms. Let
%$u_1,u_2, \ldots, u_m$ and $v_1,v_2, \ldots, v_m$
%be two concentric cycles in the plane. 

The {\em $v$-vector} $v(G)=(\ldots, v_i, \ldots)$ 
of a polyhedron $G$ enumerates
the numbers $v_i$ of vertices of degree $i$.
A polyhedron is {\em simple} if $v_i=0$ for $i\not= 3$.
The {\em $p$-vector} $p(G)=(\ldots, p_i, \ldots)$ 
of a polyhedron $G$ counts
the numbers $p_i$ of faces having $i$ sides.
For a connected plane graph $G$, we 
denote its {\em plane dual graph} by $G^*$
and define it on the set of faces of $G$ with two
faces being adjacent if they share an edge. Clearly, 
$v(G^*)=p(G)$ and $p(G^*)=v(G)$.

\bigskip

\begin{figure}
\begin{center}
\epsfxsize=90mm
\epsfbox{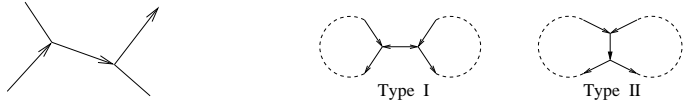}
\end{center}
\caption{A fragment of a zigzag and the two types of self-intersection.}
\label{TypeI-TypeII}
\end{figure}

%A graph $G$ is called {\em Eulerian} if there is a closed walk 
%traversing every edge of $G$ exactly once.
%It is a well-known fact that a graph is Eulerian if and only if it is
%connected and every vertex has even degree.
Let $G$ be a plane graph, such that every of its vertices has at least
degree $3$.
Then all edges, incident to a vertex $x$,
can be labeled in counter-clockwise order as
$e_1, e_2, \ldots, e_k$, where $k$ is the degree 
of $x$ in $G$. 
For any edge $e_i$, $1 \le i \le k$, the edges $e_{i+1}$,  
and $e_{i-1}$ (with $i+1$ and $i-1$ being addition modulo $k$) 
are called, respectively, the {\em left} and the {\em right}.
A circuit of edges of $G$ is
called a {\em zigzag} (or a {\em Petrie path \cite{Cox}}, 
{\em geodesic \cite{GrMo}}, {\em left-right path \cite{Sh}}), 
if, in tracing the circuit, we alternately select, as the next edge,
the left neighbor
and the right neighbor (see Figure \ref{TypeI-TypeII}).
In a $3$-valent plane graph, any pair 
of edges sharing a vertex define a zigzag.

%\begin{floatingfigure}{6cm}
%\epsfxsize=40mm
%\epsfbox{SampleZigZag.eps}
%\caption{An example of a zigzag}
%\label{fig:AzigzagSample}
%\end{floatingfigure}

Given an edge, there are two possible directions for extending
it to a zigzag. Hence, each edge is covered exactly twice
by zigzags.
A zigzag $Z$ is called {\em simple} if it has no self-intersection.
Otherwise, it has at least one edge, say, $e$, in which it 
self-intersects. Let us choose an orientation on $Z$; call the edge $e$
of {\em type I} or {\em type II} \cite{GrMo} if $Z$ traverses 
it twice in the opposite or in the same direction, respectively (see Figure 
\ref{TypeI-TypeII}). Clearly, 
the type of edge does not depend on the choice of orientation.
Edges of type I and type II correspond to edges of
{\em cocycle} and {\em cycle character} in \cite{Sh}.
The {\em signature} of a zigzag is the pair $(\alpha_1,\alpha_2)$,
where $\alpha_1$ and $\alpha_2$ are the numbers of its edges of type I and
type II, respectively.

The {\em $z$-vector} of a graph $G$ is the vector enumerating the 
lengths of all its zigzags with their signature as subscript.
The simple zigzags are put in the beginning, in increasing order 
of length, without their signature $(0,0)$, and separated by 
a semicolon from others. Self-intersecting zigzags are also ordered by
increasing lengths. If there are $m>1$ zigzags of the same
length $l$ and the same signature $(\alpha_1, \alpha_2)\not= (0,0)$, 
then we write $l_{\alpha_1,\alpha_2}^m$.
For a zigzag $Z$, its {\em intersection vector} 
$Int(Z)=(...,c_k^{m_k},...)$ is such that 
$(\ldots,c_k,\ldots)$ is an increasing sequence of sizes $c_k$ of its 
intersection with all others zigzags, and $m_k$ denote respective
multiplicities.

%This property holds also for maps on oriented surfaces.

The length of each zigzag is even, because we consecutively take
left and right turns. 
But the length of intersection
of two simple zigzags can be odd (see Klein and Dyck maps on the
Table \ref{tab:t1}; they both have oreinted genus $3$ and come as 
quotients of regular triangulations of the hyperbolic plane having
valency $7$ and $8$, respectively). Returning to our plane case, the
number of intersection can be even only.

%The length of each zigzag is even, because each pair of zigzags intersect 
%in even number of edges and any edge of self-intersection is counted twice 
%in the length of zigzag. If, moreover, the graph $G$ has only simple zigzags, 
%then the number of edges of $G$ is even; for a $3$-valent graph $G$ it implies
%that its number of vertices is divided by $4$. This evenness property does
%not hold for maps on oriented surfaces; for example, $z$-vector
%of the regular Klein map (having oriented genus $3$, valency $7$, $24$ vertices
%and $56$ triangular faces) is $21^8$ (NEED TO CHECK) and the intersection 
%vector of each zigzag is $3^7$.

A graph is called {\em $z$-uniform} if all zigzags have the same
length and signature, {\em $z$-transitive} if
its group acts transitively on zigzags. Clearly, 
$z$-transitivity implies $z$-uniformity. A graph is called {\em $z$-knotted}
if it has only one zigzag, {\em $z$-balanced} if all its zigzags of the
same length and same signature, have identical intersection vectors.
Clearly, $z$-transitivity implies $z$-balancedness.

Let $Z_1$, \dots, $Z_p$ be all zigzags of $G$. On every zigzag
choose an orientation; there are $2^p$ 
possibilities. Every edge, being covered twice by zigzags (one or two), 
can be, with respect to the choosen orientation of all zigzags, oriented 
in the opposite or same direction. We call such an edge of {\em type I}
or {\em type II with respect to the fixed orientation}. If an edge
is a self-intersection edge of some zigzag, then the notion of type
is independent of the fixed orientation and coincides with the notion 
of type I and II, which was introduced above.

While for dual graphs $v$- and $p$-vectors are interchanged,
the $z$-vector remains the same, except that type I
and type II in the signature are interchanged.
Other kinds of dualities (for example, interchanging
$z$- and $p$-vector, but preserving $v$-vector) are considered 
in \cite{L}.
%NEED REFERENCE FOR INVARIANCE OF ZIGZAG BY DUALITY

\bigskip

The {\em medial} graph $Med(G)$ of $G$ is defined 
on the set of edges of $G$ with two edges being adjacent
if they share a common vertex and if they belong
to a common face. The medial graph is $4$-valent and 
$Med(G)=Med(G^*)$.
Zigzags of $G$ and $G^*$ correspond to central 
circuits (see below) of $Med(G)$.
Any $4$-valent plane graph $H$ is the medial graph
for a pair of mutually dual plane graphs: one can assign two 
colors to the faces of $H$ in the ``chess way'', such that
 no two adjacent faces of $H$ have the same color. Two faces
 of the same color are said to be adjacent if they share
 a vertex.

In general, a {\em central circuit} of
an {\em Eulerian} (i.e., having only vertices of even 
degree) plane graph is a circuit, which is obtained
by starting with an edge and continuing at each vertex by the 
edge opposite to the entering one. A {\em link} is one or more 
circles (components of a link) embedded into the space $\R^3$;
a link with one component is called {\em knot}.
A {\em projection} of a link is a drawing of it on the plane with gaps
representing underpass and solid line representing
overpass. If a drawing of a link has alternating
underpasses and overpasses, then it is called {\em alternating};
so, we can see it just as a $4$-valent plane graph.
We will consider only {\em minimal} projections, i.e., those
without $1$-gons (loops).
Clearly, each edge belongs to exactly one central circuit and any 
$4$-valent graph without $1$-gons can be seen as a
minimal projection of an alternating link with components 
corresponding to its central circuits. 
Since zigzags in $G$ correspond to central circuits 
in $Med(G)$, each zigzag corresponds to a component
of the corresponding alternating link.
The bipartition of edges of a $z$-knotted graph into type I
and type II corresponds to a bipartition of the vertices of its
medial graph (which is a projection of an alternating knot).

Call a polyhedron {\em two-faced} if it has only $a$- and $b$-gonal 
faces with $3\leq a < b\leq 6$. In particular, denote $3$-valent
two-faced $n$-vertex polyhedron by $3_n$, $4_n$, $5_n$ if $(a,b)$=$(3,6)$,
$(4,6)$, $(5,6)$, respectively.
Given two circuits $u_1, \dots, u_m$ and $v_1, \dots, v_m$,
an {\em $m$-sided prism} $Prism_m$ is formed when every $u_i$ is joined to
$v_i$ by an edge. Now, an {\em $m$-sided antiprism}
$APrism_m$ is formed by adding the cycle $u_1$,$v_2$,$u_2$,$v_3$,
\ldots, $v_m$,$u_m$,$v_1$,$u_1$.
%The Conway graph $(k\times m)^*$ (see, for example, \cite{Kaw}) is, 
%for $k=2$, $m$-sided antiprism; for $k>2$, it comes from 
%$((k-1)\times m)^*$ by inscribing an $m$-gon in the first of its
%two $m$-gons.

A {\em point group} is a finite subgroup of the group $O(3)$ of isometries of
the space $\R^3$, fixing the origin. The connection with plane graphs come
from representing them on the sphere centered at the origin.
The list of point groups is splitted into two classes: the infinite families
and the sporadic cases. Every point group $G$ contains a normal subgroup
formed by the rotations of $G$.
The group, denoted $C_m$, is the cyclic group of rotations by angle
$\frac{2\pi}{m}k$ with $0\leq k\leq m-1$ around a fixed axis $\Delta$.
Both groups $C_{mv}$ and $C_{mh}$ contains $C_{m}$ as normal subgroups
of rotations. The group $C_{mh}$ (respectively, $C_{mv}$) is the group,
generated by $C_{m}$ and a symmetry of plane $P$, with $P$ being orthogonal
to $\Delta$ (respectively, containing $\Delta$).
The group $D_m$ is the group, generated by $C_m$ and a rotation by angle $\pi$,
whose axis is orthogonal to $\Delta$. The point group $D_{mh}$ (respectively,
$D_{md}$) is generated by $C_{mv}$ and a rotation by angle $\pi$, whose axis
is orthogonal to $\Delta$ and contained in a plane of symmetry (respectively, 
going between two planes of symmetry). Both, $D_{mh}$ and $D_{md}$, contain
$D_{m}$ as a normal subgroup. If $N$ is even, one defines the point group $S_N$
to be the cyclic group generated by the composition of a rotation by angle
$\frac{2\pi}{N}$ with axis $\Delta$ and a symmetry of plane $P$
with $P$ being orthogonal to $\Delta$. The particular cases $C_1$, $C_{s}$
and $C_{i}$ correspond to the trivial group, the plane symmetry group
and the central symmetry inversion group.
The point groups $T_d$, $O_h$ and $I_h$ are the symmetry group of the
Tetrahedron, Cube and Icosahedron; the point groups $T$, $O$ and $I$
are their normal subgroup of rotations. The point group $T_{h}$ is
formed by all $f$ and $-f$ with $f\in T$. More detailed description
of point groups are available, for example, from \cite{FM} and \cite{D1}.

\begin{table}
\scriptsize
\begin{center}
\begin{tabular}{||c|c|c|c||}
\hline
\hline
\#\,edges&Polyhedron& $z$-vector&Int. vector \\ \hline
\hline
6  &{\em Tetrahedron}&$4^3$&$2^2$\\ 
12 &{\em Cube}&$6^4$&$2^3$\\ 
12 &{\em Octahedron}&$6^4$&$2^3$\\ 
30 &{\em Dodecahedron}&$10^6$&$2^5$\\ 
30 &{\em Icosahedron}&$10^6$&$2^5$\\ \hline
24 &{\em Cuboctahedron}&$8^6$&$2^4,0$ \\%=Med(cube)=Med(octahedron)
60 &{\em Icosidodecahedron}&$10^{12}$&$2^5,0^6$ \\%=Med(ico.)=Med(dode.)
48 &{\em Rhombicuboctahedron}&$12^8$&$2^6,0$\\ %=Med(cuboctahedron)
120&{\em Rhombicosidodecahedron}&$20^{12}$&$2^{10},0$\\ %=Med(icosidodecahedron)
72 &{\em Truncated Cuboctahedron}&$18^8$&$6,2^6$\\ 
180&{\em Truncated Icosidodecahedron}&$30^{12}$&$10,2^{10}$\\ 
18 &{\em Truncated Tetrahedron}&$12^3$&$6^2$\\ 
36 &{\em Truncated Octahedron}&$12^6$&$4,2^4$\\ 
36 &{\em Truncated Cube}&$18^4$&$6^3$\\
90 &{\em Truncated Icosahedron}&$18^{10}$&$2^9$\\ 
90 &{\em Truncated Dodecahedron}&$30^6$&$6^5$\\ 
60 &{\em Snub Cube}&$30_{3,0}^4$&$8^3$\\ 
150&{\em Snub Dodecahedron}&$50_{5,0}^6$&$8^5$\\ \hline
3m &$Prism_m$, $m\equiv 0 \pmod{4}$&$(\frac{3m}{2})^4$&$(\frac{m}{2})^3$\\
3m &$Prism_m$, $m\equiv 2 \pmod{4}$&$({3m}_{\frac{m}{2},0})^2$&$2m$\\ 
3m &$Prism_m$, $m\equiv 1,3 \pmod{4}$&${6m}_{m,2m}$&\\ 
4m &$APrism_m, m \equiv 0 \pmod{3}$&$(2m)^4$&$(\frac{2m}{3})^3$\\ 
4m &$APrism_m, m \equiv 1,2 \pmod{3}$&$2m;6m_{0,2m}$&\\ \hline\hline
84 &{\em Klein map}                  &$8^{21}$&$1^8, {0}^{12}$\\
48 &{\em Dyck map}                   &$6^{16}$&$1^6, {0}^9$\\
\hline
\hline
\end{tabular}
\end{center}
\caption{Zigzag structure of Platonic and semiregular polyhedra; also, of Klein and Dyck maps.}
\label{tab:t1}
\end{table}

In Table \ref{tab:t1}, zigzag notions are illustrated by the
regular and {\em semiregular} polyhedra (i.e., such that their symmetry group
is transitive on vertices).
This table is a slightly more precise version (indicating types
of self-intersection of zigzags) of Table 1 of \cite{DHL}.
It is easy to see that for self-intersecting zigzags of prisms and
antiprisms, the self-intersections of type I, if any, occur on the
edges of the two $m$-gons, while those of type II occur on other
edges.

The notion of a zigzag was generalized to locally-finite 
{\em infinite} planar $3$-connected graphs. All such 
{\em edge-transitive} graphs without self-intersecting zigzags 
(i.e., simple circuits and doubly infinite rays are the only 
zigzags) were classified in \cite{GrSh87}: these are the three regular 
partitions $(6^3)$, $(3^6)$, $(4^4)$ of the Euclidean plane, 
the Archimedean 
partition $(3.6.3.6)$, its dual and several infinities 
of partitions of the hyperbolic plane. Also the notion 
of zigzag (Petrie polygon) was extended in 
\cite{Cox} p.~223 for $n$-dimensional polytopes and honeycombs.

\bigskip

This paper is based on extensive computation; in particular:
\begin{itemize}
\item[(i)] All computations used the GAP computer algebra system
  \cite{Gap} and the package PlanGraph (\cite{D2}) by the second author.
\item[(ii)] The program CPF (\cite{TH}) by Harmuth was used 
to generate the graphs of type $3_n$ and $4_n$. In fact, the name $n_{x}$ 
for a graph indicates that a graph appears at $x$-th position in the 
output of CPF. Also, CaGe (\cite{CaGe}) was used for graph drawings.
\item[(iii)] For graphs of type $5_n$, we used also the
face-spiral algorithm notation given in \cite{FM}; a computation by
Brinkmann (\cite{B}) was used in Table \ref{tab:t3}.
\end{itemize}

\section{General results for plane graphs}

\begin{theorem}\label{Shtogrin-Shank}
For any planar bipartite graph $G$ there exist an orientation of zigzags, with respect to which each edge has type I.
\end{theorem}
\proof We represent the graph $G$ on the sphere. The list of vertices, adjacent to a given vertex, can be organized in counter-clockwise order.

Let the vertex-set $V$ be partitioned into the two subsets $V_1$ and
$V_2$ of the bipartition. Fix a
zigzag $Z$; it will turn left at vertices, say, $v\in V_1$ and
right at vertices $v\in V_2$. It is easy to see that the edges of
self-intersection of $Z$ can be only of type I.

Take another zigzag $Z'$, having a common edge $e$ with $Z$. We choose an orientation on $Z'$, such that $e$ is an edge of type I. Then $Z'$ will turn left at vertices $v\in V_1$ and right at vertices $v\in V_2$. Iterating this construction, all zigzags will be oriented and all edges will have type I with respect to this orientation of zigzags. \qed

\vspace{0.1cm}
In the case of one zigzag, this result was already known (\cite{Sh}). 
This theorem is also valid for any bipartite graph, which is embedded
in an oriented surface, in view of the well-known topological fact that any
two-dimensional orientable manifold admits coherent orientation of its
faces.

\begin{proposition}
(\cite{DDF}) Let $G$ be a plane graph; for any orientation of all zigzags of $G$, we have:

(i) The number of edges of type II, which are incident to any fixed {\em vertex}, is even.

(ii) The number of edges of type I, which are incident to any fixed {\em face}, is even.
\end{proposition}
\proof (i) For each vertex the number of times, that a zigzag enters it, should be equal to the number of times that a zigzag leaves it.

(ii) Passing to the dual graph, the type of edges are interchanged. \qed

A consequence of the above proposition for any $3$-valent plane graph $G$ (see
\cite{DDF}) is that the set of vertices can be partitioned into two
classes, say, class I and class II, where class I stands for vertices,
which are incident to three edges of type I, and class II stands for
vertices, which are incident to one edge of type I and two edges of
type II. If $n$ denotes the number of vertices of $G$ and $n_1$
denotes the number of its vertices of class I, then the number of edges of
type I is $n_1+\frac{n}{2}$. Clearly, the number of edges of type II
is equal to the number of vertices of class II, i.e., to $n-n_1$.

It was conjectured in (\cite{DDF}) that any $3$-valent plane graph,
which is $z$-knotted, has an odd number of edges of type I. The condition 
of $3$-valency is necessary: for example, the
$3$-bipyramid (i.e., $Prism_3^*$) has $z$-vector $18_{6,3}$.

\bigskip

\begin{figure}
\centering
\begin{minipage}[b]{4.5cm}
\centering
\epsfig{figure=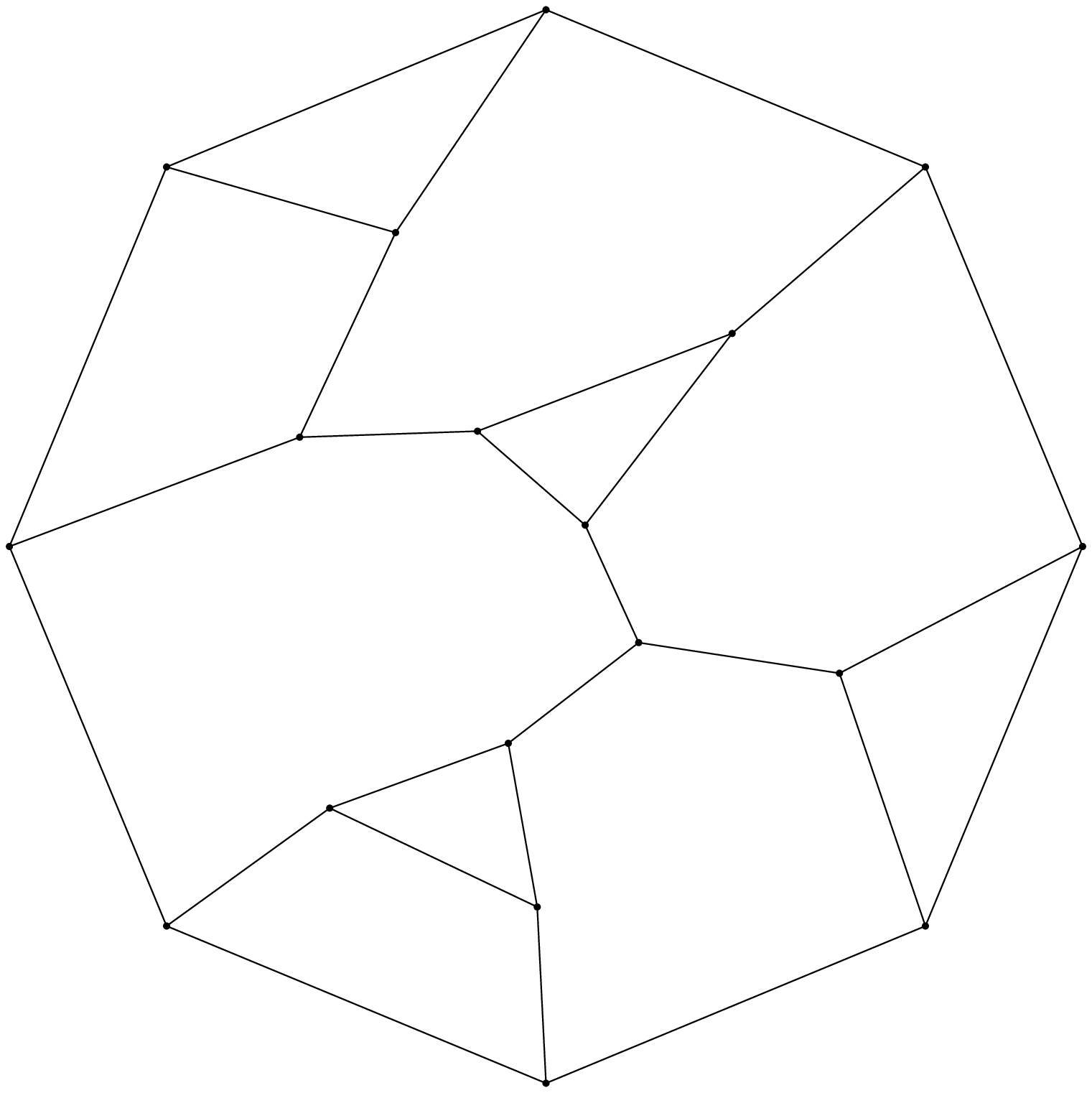,height=4cm}\par
$18$-vertices graph $(C_2)$, $z=8^3; 30_{2,5}$
\end{minipage}
\begin{minipage}[b]{4.5cm}
\centering
\epsfig{figure=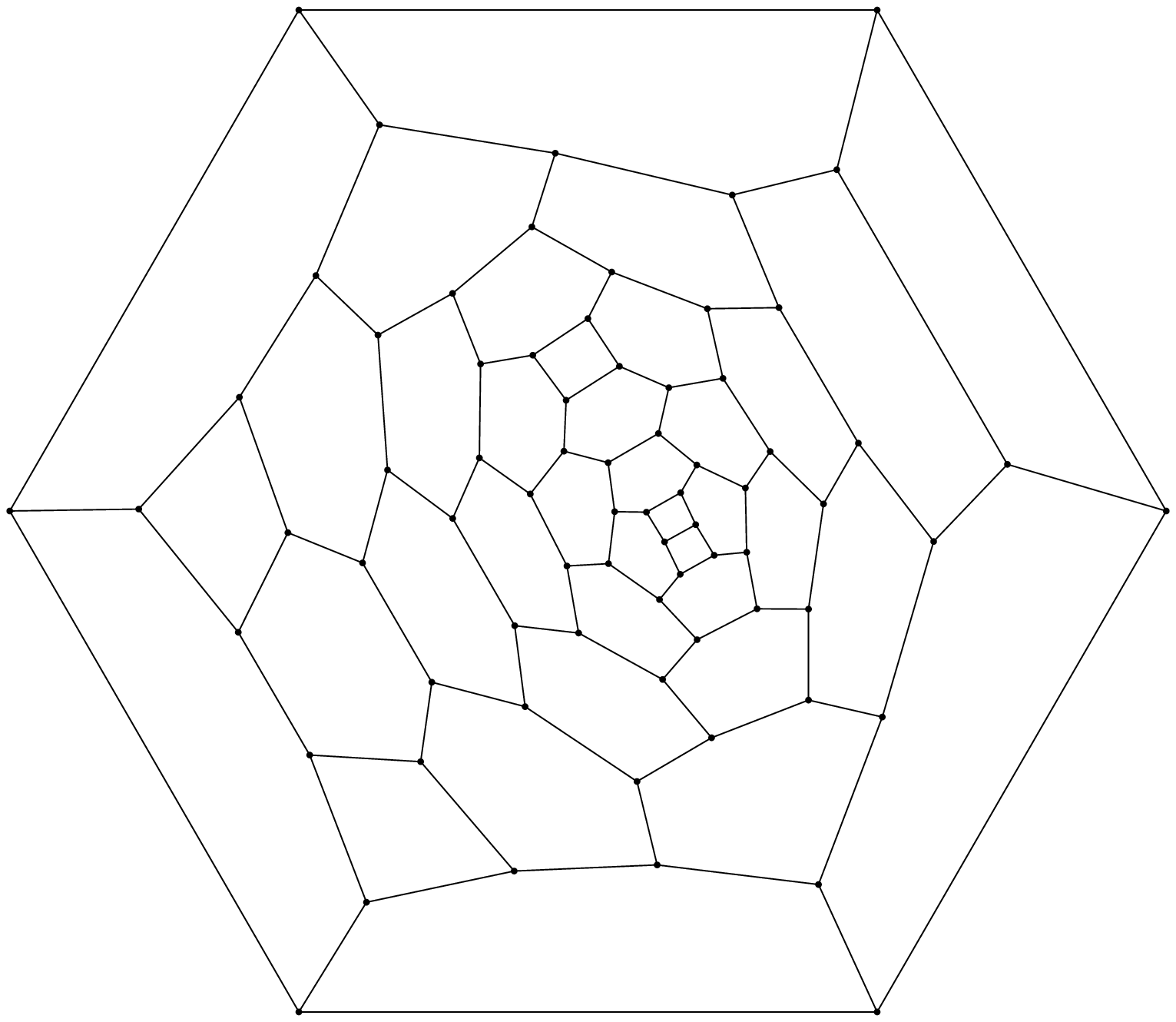,height=4cm}
Graph $4_{72}$($C_1$), $z=30_{1,0}, 54^2_{4,0}, 78_{13, 0}$
\end{minipage}
\begin{minipage}[b]{4.5cm}
\centering
\epsfig{figure=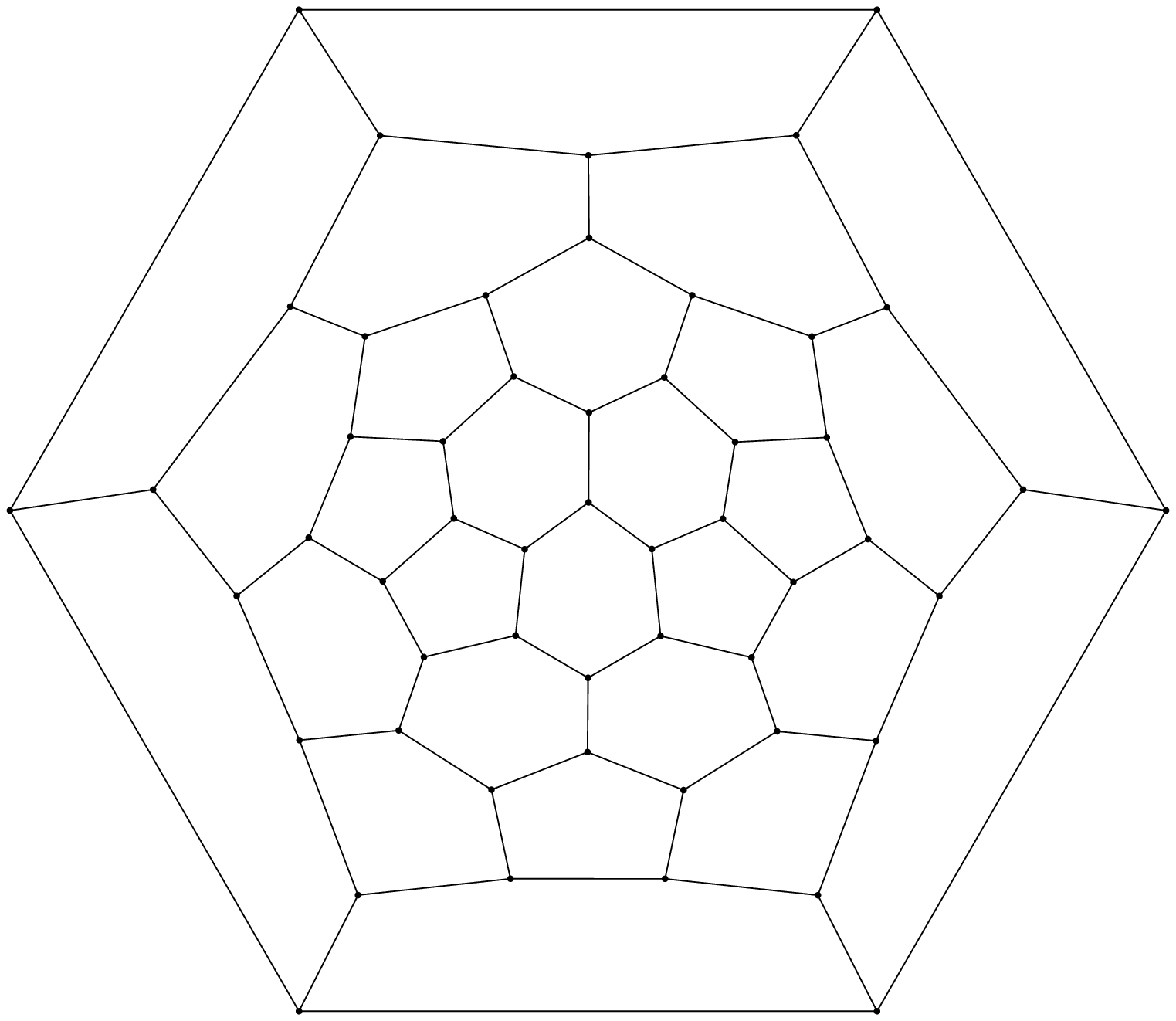,height=4cm}
Graph $5_{52}$($D_{2d}$), $z=16^4; 92_{12, 12}$
\end{minipage}
\caption{Smallest $z$-unbalanced $3$-valent graphs.}
\label{fig:Unbalanced_graphs}
\end{figure}

Clearly, any graph with $1$, $2$ or $3$ zigzags is $z$-balanced.
See smallest $z$-unbalanced $3$-valent graph and graphs of type $4_n$, $5_n$ 
on Figure \ref{fig:Unbalanced_graphs}. We did computations trying 
to find an example of a $z$-unbalanced and {\em $z$-uniform} $3$-valent 
polyhedron (for $n\leq 22$), graph of type $4_n$ (for $n\leq 250$) and 
$5_n$ (for $n\leq 80$); surprisingly, we did not find any
such example.

%\begin{figure}
%\centering
%\mbox{\subfigure[{$18$-vertices $3$-valent graph $(C_2)$}]{\epsfig{figure=FirstUnbalancedTrig.ps,width=.27\textwidth}}\quad
%\subfigure[{Graph $4_{72}$($C_1$)}]{\epsfig{figure=FirstUnbalanced4n.ps,width=.30\textwidth}}\quad
%\subfigure[{Graph $5_{52}$($D_{2d}$)}]{\epsfig{figure=FirstUnbalanced5n.ps,width=.30\textwidth}}}\caption{Smallest
%unbalanced $3$-valent graphs}
%\label{fig:Unbalanced_graphs}
%\end{figure}

%The first $3$-valent plane $3$-connected
%graph, which is not $z$-balanced, is
%\begin{center}
%\epsfxsize=60mm
%\epsfbox{
%\end{center}
%with $z=8\sp{3}; 30_{2,5}$ and whose intersection vectors of $3$ zigzags of length $8$ are $(4,6^2)$, $(0,2,6)$ and $(2^2,4)$.

In the rest of this section, we give a local Euler formula for zigzags.
Let $G$ be a plane $3$-valent graph. Consider a patch $A$ in $G$, 
which is bounded by $t$ arcs (i.e., paths of edges) belonging to 
zigzags (different or coinciding). 
%Call a {\em vertex} of the
%patch any edge of intersection of two consecutive arcs of the 
%boundary.

\begin{figure}
\centering
\epsfxsize=100mm
\epsfbox{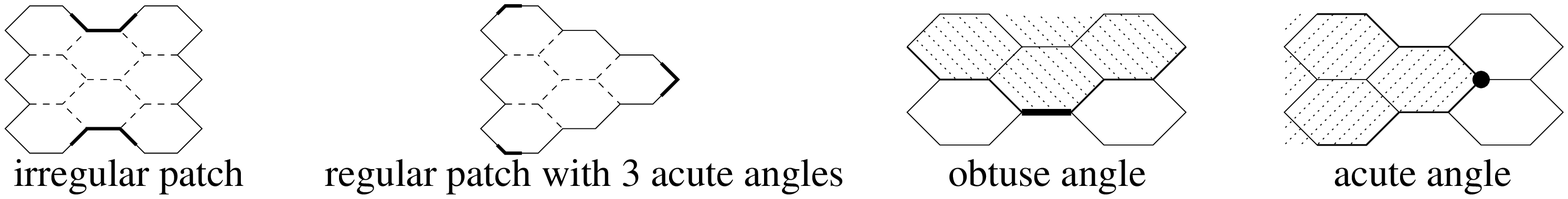}
\caption{Example of patches and their angles.}
\label{fig:Regular_unregular}
\end{figure}

We admit also $0$-gonal patch $A$, i.e., just the interior of a simple zigzag. 
Suppose that the patch $A$ is {\em regular}, i.e., the continuation 
of any of its bounding arcs (on the 
zigzag, to which it belongs) lies outside of the patch (see Figure 
\ref{fig:Regular_unregular}). Let 
$p'(A)=(...,p'_i,...)$ be the $p$-vector enumerating the faces
of $G$, which are contained in the patch $A$.

There are two types of intersections of arcs on the boundary of a regular 
patch: either intersection in an edge of the boundary, or intersection in 
a vertex of the boundary. Let us call these types of intersections
{\em obtuse} and {\em acute}, respectively (see Figure
\ref{fig:Regular_unregular}); denote by $t_{ob}$ and $t_{ac}$ the 
respective number of obtuse and acute intersections. Clearly, 
$t_{ob}+t_{ac}=t$, where $t$ is the number of arcs forming the patch.

%\begin{center}
%\epsfxsize=60mm
%\epsfbox{LocalDrawing.eps}
%\end{center}

\begin{theorem}\label{Local-Euler-Formula}
Let $G$ be a $3$-valent plane graph ($1$-gonal and $2$-gonal faces are
permitted). If $A$ is a regular patch, then the following equality holds
$$6-t_{ob}-2t_{ac}=\sum_{i\geq 1} (6-i)p'_{i}\; .$$

\end{theorem}
\proof Let $P(A)$ be the induced plane subgraph of $G$ formed by the
patch $A$. The vertices of $P(A)$ have degree $2$ or $3$; we
denote their respective number by $v_2$ and $v_3$. The exterior face is
$l$-gonal, where $l$ is the length of the boundary of the patch.
Direct enumeration gives the following expressions for the number $|E|$ of edges of $P(A)$
$$|E|=\frac{1}{2}(3v_3+2v_2)=\frac{1}{2}(l+\sum_{i\geq 1}ip'_{i})\;.$$

Euler's formula, applied to the plane graph $P(A)$, yields
$$2=(v_2+v_3)-|E|+(1+\sum_{i\geq 1}p'_{i})\;.$$

The above two expressions of $|E|$ give
$$1=-\frac{v_3}{2}+\sum_{i\geq 1}p'_{i}=(v_2+v_3)-\frac{l}{2}+\sum_{i\geq 1} (1-\frac{i}{2}) p'_{i}\;.$$

Eliminating $v_3$ yields
$3=v_2-\frac{l}{2}+\sum_{i\geq 1} (3-\frac{i}{2}) p'_{i}\;.$

For example, if the patch $A$ is a $0$-gon (i.e., if it is bounded by a simple zigzag), then $l=2v_2$ and we get 
$$3=\sum_{i\geq 1} (3-\frac{i}{2}) p'_{i}.$$

Denote by $v'_2$ and $v'_3$ the numbers of vertices on the boundary of $A$, having degree $2$ and $3$, respectively. Clearly, $l=v'_2+v'_3$ and $v'_2=v_2$.

Now we use that $A$ is a regular patch. Specifying acute and obtuse
types of arc intersections on the boundary of $A$, we get $v'_2=v'_3+t_{ob}+2t_{ac}$, which gives
$6-t_{ob}-2t_{ac}=\sum_{i\geq 1} (6-i) p'_{i}\;.$ \qed

\begin{corollary}\label{Easy-consequence-Euler-Formula}
Let $G$ be a $3$-valent plane graph. Let $A$ be a $t$-gonal regular patch with $p$-vector $(\dots, p'_{i}, \dots)$ in $G$. Then we have:

(i) If $G$ has no $q$-gonal faces with $q>6$, then $t\leq 6$.

(ii) If $G$ is a two-faced graph $w_{n}$, $w\in \{ 3,4,5\}$, then
$$6-t_{ob}-2t_{ac}=(6-w)p'_w\;.$$

\end{corollary}
\proof (i) By the above theorem 
$$6-t_{ob}-2t_{ac}=\sum_{1\leq i\leq 6}(6-i)p'_i\geq 0\;,$$
which gives the result.
(ii) also follows by direct application of Theorem \ref{Local-Euler-Formula}. \qed

\section{Railroads and pseudo-roads}
Let $G$ be a plane graph. Call a {\em railroad} of $G$ a 
circuit of hexagonal faces in $G$, 
such that each hexagon is adjacent to its neighbors on opposite edges
(in \cite{GrMo}, a {\em simple}, i.e., without self-intersections,
railroad is called a {\em belt}). Clearly, 
a railroad is bounded by two ``parallel'' (concentric if the railroad is simple) zigzags.
The graph $G$ is called {\em tight} if it has no railroads. 
A railroad in $G$ corresponds to a central circuit in $G^*$, 
which goes only through $6$-valent vertices.

We associate to each railroad a {\em representing} plane 
curve in the following way: in each of its hexagons one
connects, by an arc, the midpoints of opposite edges, on which
it is adjacent to its two neighbors. The sequence of those arcs
can be seen as a Jordan curve in the 
plane and self-intersections of railroad correspond to self-intersections 
of the curve.

Those self-intersections can be only double or triple, because 
a railroad consists of hexagons. If the railroad has no triple 
self-intersection points, then this curve can be seen as a projection 
(not minimal, if there are $1$-gons) of an alternating knot 
(see, for example, \cite{Kaw} and \cite{Rolf}) with $n$ crossings, 
where $n$ is the number of self-intersection points. If there are
several railroads without triple intersections and triple 
self-intersections, then the set of plane curves, representing them,
can be seen as a projection of an alternating link.
%A link diagram is called {\em alternating} if, after the choice 
%of an orientation, an over-crossing and an under-crossing 
%appear alternately as the arcs of the link are traversed.

A set of railroads can be seen as a plane graph $H$
with valency $4$ or $6$ of its vertices, accordingly to double or triple 
points of intersection or self-intersection of the curves. Every
$t$-gonal face $F$ of $H$ can be viewed as a regular patch 
with $t$ intersections of arcs.
In the case, when $G$ has no $q$-gonal faces with $q>6$, Corollary 
\ref{Easy-consequence-Euler-Formula}.(i) implies that $t\leq 6$.

\begin{figure}
\centering
\epsfxsize=100mm
\epsfbox{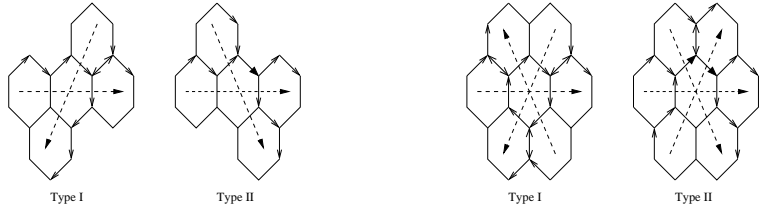}
\caption{Types of self-intersection of a railroad.}
\label{fig:TypeSelfIntersection}
\end{figure}

A railroad $R$ is bounded by two zigzags $Z_1$ and $Z_2$, which have the
same length and signature. Each edge of self-intersection of $Z_1$, $Z_2$
corresponds to a self-intersection of $R$ and so, to an hexagon of
self-intersection of $R$. Since $Z_1$ can self-intersect in an edge
of type I or type II, there are two types of double self-intersections
of railroad. Simple analysis yields two types for triple 
self-intersections of railroad: either $Z_1$ self-intersects in three
edges of type I or it self-intersects in one edge of type I and two 
edges of type II. See Figure \ref{fig:TypeSelfIntersection} for all 
possibilities.

\begin{figure}
\centering
\begin{minipage}[b]{4.6cm}%
\centering
\epsfig{figure=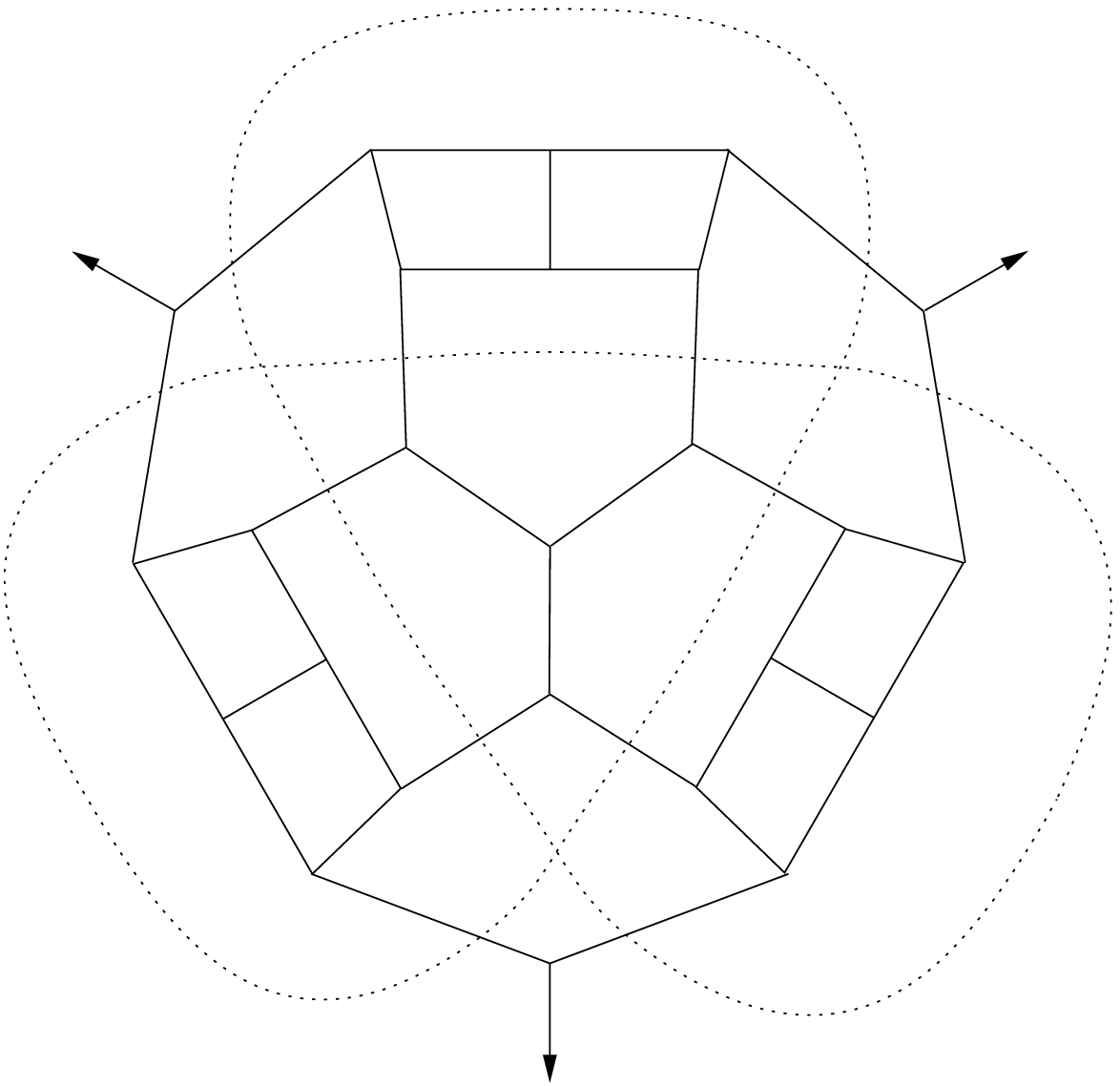,height=4cm}
Graph $26_{2}(D_{3h})$
\end{minipage}
\begin{minipage}[b]{4.6cm}%
\centering
\epsfig{figure=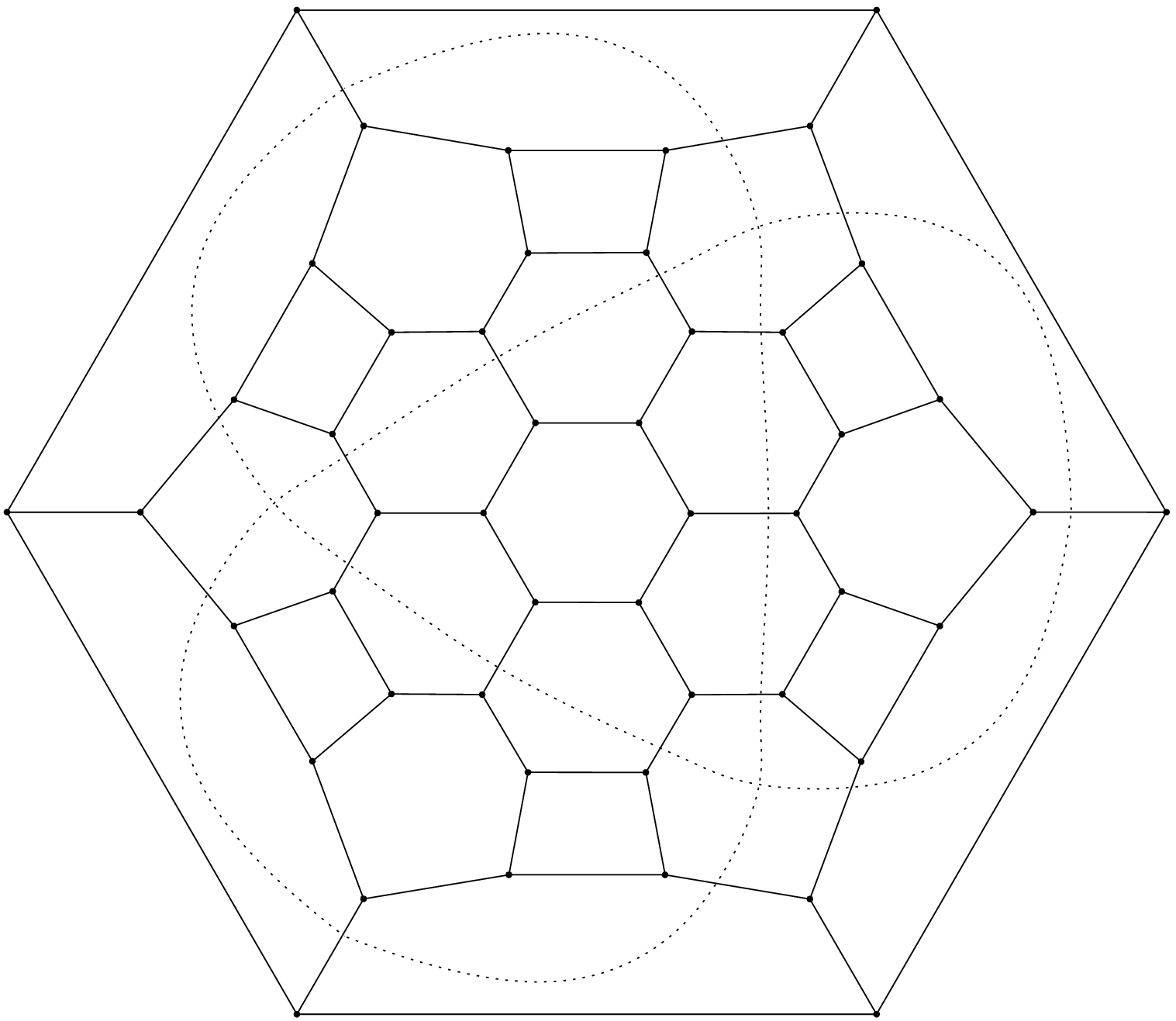,height=4cm}
Graph $48_{10}(D_{6h})$
\end{minipage}
\begin{minipage}[b]{4.6cm}%
\centering
\epsfig{figure=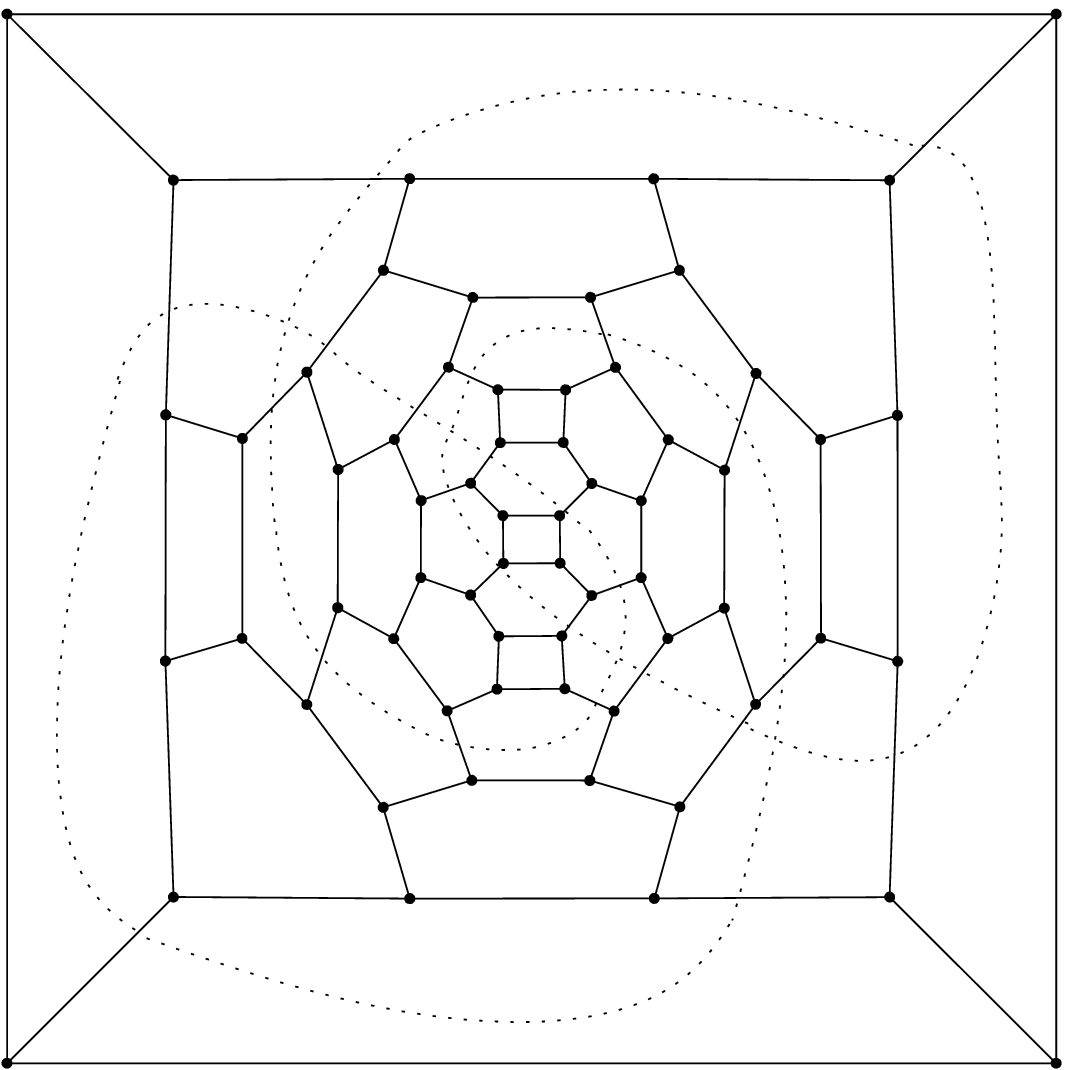,height=4cm}
Graph $64_{17}(D_{2d})$
\end{minipage}
\caption{Three graphs of type $4_n$ with self-intersecting railroads.}
\end{figure}

\begin{table}
\scriptsize
\begin{equation*}
\begin{array}{||c|c|c|c|c|l||}
\hline\hline
Knot       &First\;\; appearance&m &p_1, \dots,p_6&Group     &z-vector\\\hline
0_1        &20_3                &0 &              &D_{3d}    &6^3; 42_{12,0}\\
3_1        &26_2                &3 &0,3,2,0,0,0   &D_{3h}    &24_{3,0}^{2}, 30_{3,0}\\
4_1        &64_{17}(twice)      &4 &0,2,4,0,0,0   &D_{2d}    &48_{4,0}^4\\
5_2        &74_{24}             &5 &0,3,2,2,0,0   &C_{2v}    &18^{4}; 48_{5,0}^{2}, 54_{3,0}\\
6_3        &74_{21}             &6 &0,2,4,2,0,0   &C_{2}     &54_{6,0}^{2}, 114_{27,0}\\
7_4        &80_{30}(twice)      &7 &0,4,2,2,0,1   &D_{2h}    &60_{7,0}^{4}\\
9_{23}     &62_{16}             &9 &0,2,6,2,0,1   &C_{2v}    &60_{9,0}^2, 66_{9,0}\\
9_{40}     &56_{19}             &9 &0,8,3,0,0,0   &D_{3h}    &12^4; 60_{9,0}^2\\
11_{332}   &88_{44}             &11&0,2,4,7,0,0   &C_{2v}    &24^{4}; 84_{11,0}^{2}\\
15_{y}     &80_{29}             &15&0,0,8,6,0,0   &D_{3d}    &12^{6}; 84_{12,0}^{2}\\
18_{y}     &80_{28}             &18&0,0,14,0,6,0  &D_{3}     &24^{2};96_{18,0}^{2}\\
\hline\hline
\end{array}
\end{equation*}
\caption{Alternating knots in railroads of graphs of type $4_n$, $n\leq 88$.}
\label{tab:List-of-knots}
\end{table}

In order to illustrate the notion of railroads, we give in the remainder 
of this section two lists: those having, or, respectively, not having, triple
self-intersections
of the curves representing railroads in two-faced $n$-vertex polyhedra
of type $4_n$ for small $n$. In Tables \ref{tab:List-of-knots} and
\ref{tab:Curve-with-triple-points}, we present the first cases of such curves;
the first column being the name of the curve.
The smallest graph of type $4_n$, which is denoted by $n_x$, $x$ being the number of the graph in CPF output (see \cite{TH}), is listed in the column
``First appearance''. In some cases the curve appears twice
in this graph; we mark this situation by putting ``twice'' in
parenthesis. 
The third and fourth column give the number $m$ of intersection points
and the number $p_i$ of $i$-gons, $1\leq i\leq 6$. 
Columns ``Group'' and ``$z$-vector'' give the point group and the $z$-vector of
the corresponding graph of type $4_n$.

The computation of all railroads in the
class of polyhedra $4_n$, $n\leq 88$, with at least 
two self-intersections, all being double self-intersections, produced a list of
projections of alternating knots presented in Table \ref{tab:List-of-knots}.
We use Rolfsen's notation, see \cite{Rolf}, for knots with at most $10$ 
crossings and those of Thistlethwaite (\cite{Thi}), for other knots.
We denote by $15_y$ and $18_{y}$ the $15$- and $18$-crossing alternating
knots, appearing in a graph of type $4_{80}$, which are $80_{29}$ and 
$80_{28}$, respectively. (The knot $15_{y}$ comes from $APrism_3$ or 
from knot $9_{40}$ by inscribing consecutively $3$ or $2$ triangles, respectively.)
All curves, representing railroads with triple self-intersections in the
graphs $4_n$, with $n\leq 142$ are given in Table \ref{tab:Curve-with-triple-points}.
The notation $i-j$, given in first column, means that $i$ is the number of triple points
of the curve and $j$ is order of appearance (amongst curves with $i$ triple points)
in this table.

\begin{table}
\centering
\begin{minipage}[b]{4.2cm}%
\centering
\epsfig{figure=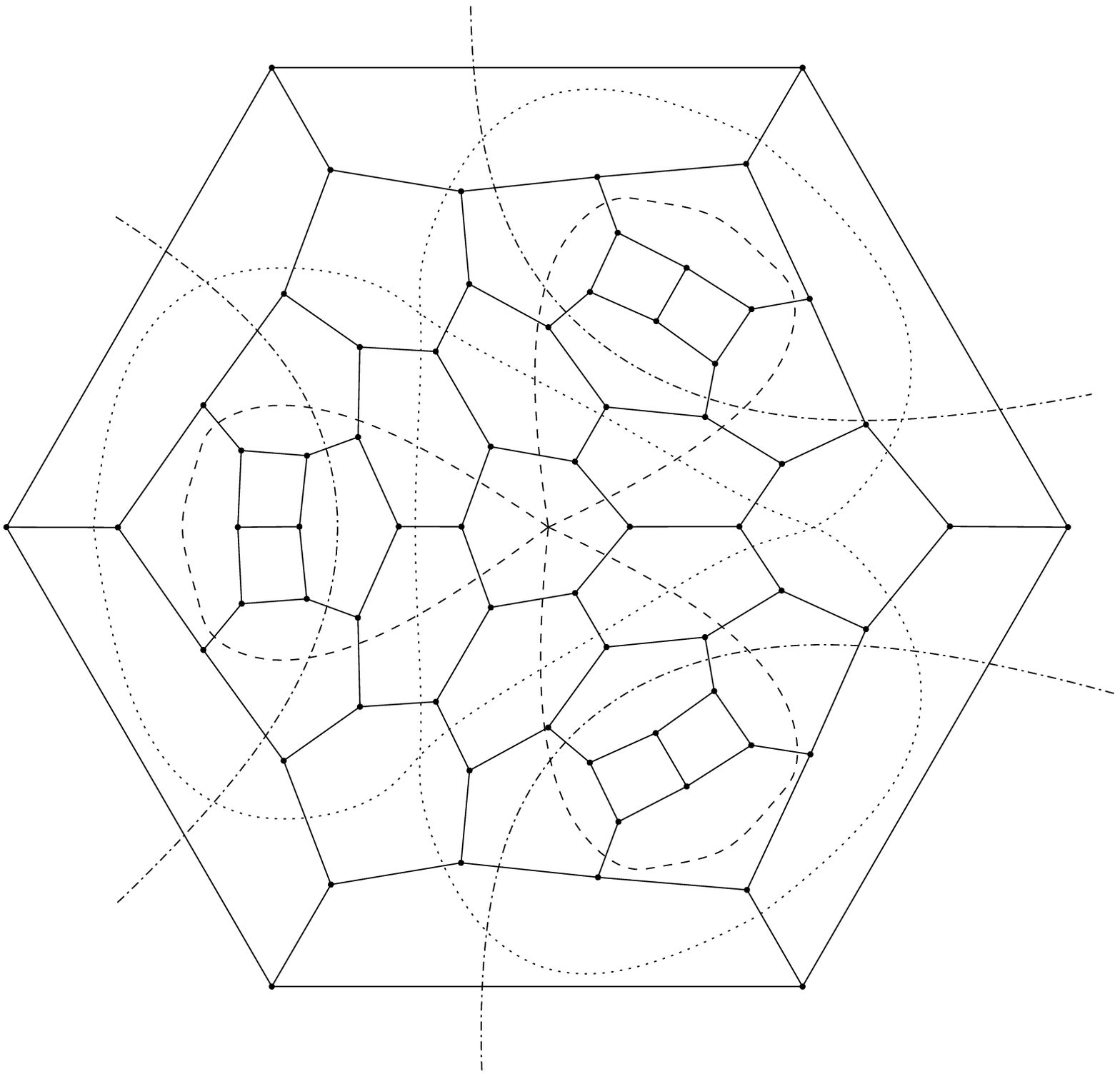,height=4cm}
Graph $66_{11}$, $1-1$
\end{minipage}
\begin{minipage}[b]{4.2cm}%
\centering
\epsfig{figure=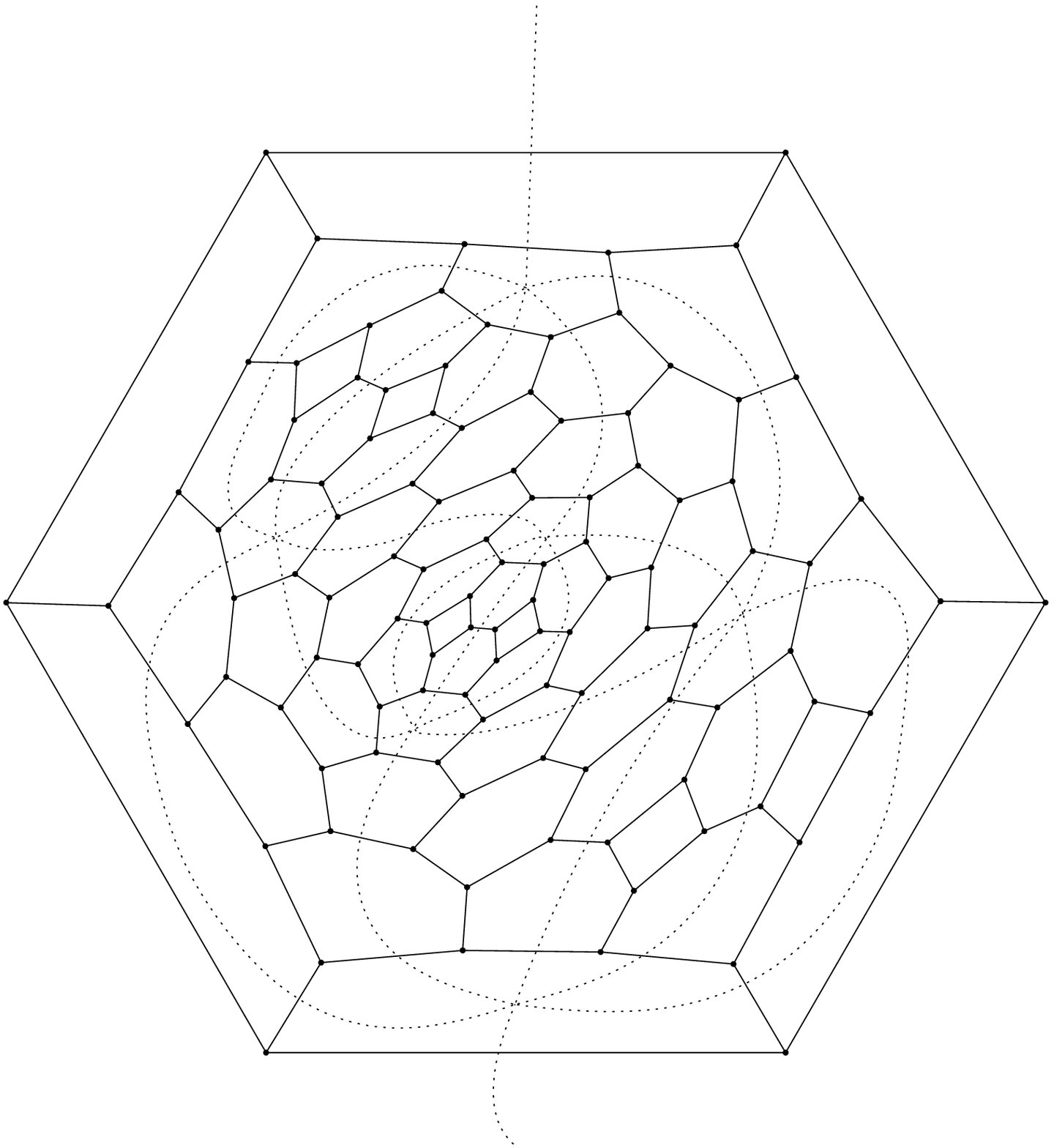,height=4cm}
Graph $108_{18}$, $6-1$
\end{minipage}
\begin{minipage}[b]{4.2cm}%
\centering
\epsfig{figure=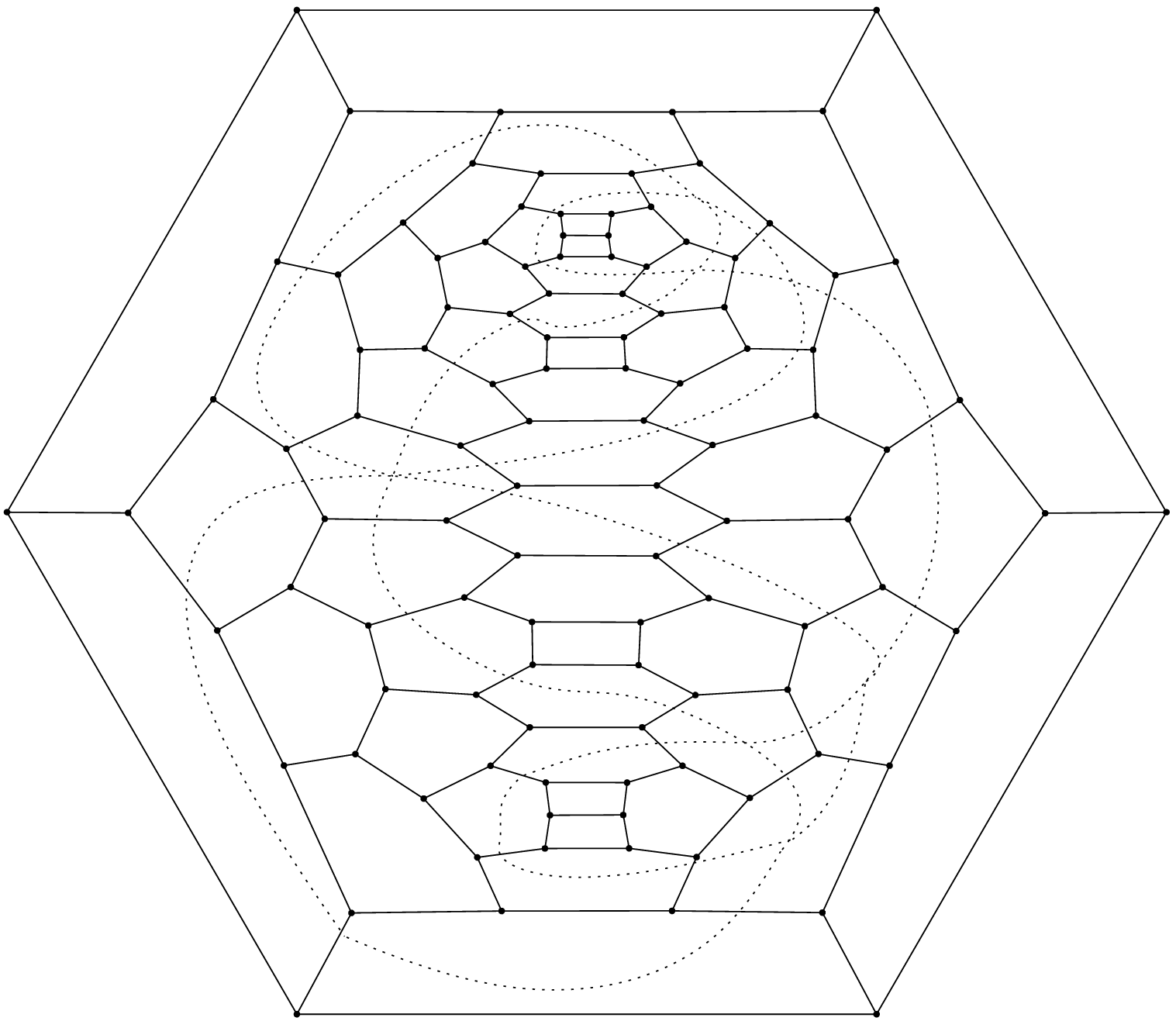,height=4cm}
Graph $108_{48}$, $1-6$
\end{minipage}
\begin{center}
Three graphs of type $4_n$ with triple self-intersecting railroads.
\end{center}
{\scriptsize
\begin{equation*}
\begin{array}{||c|c|c|c|c|l||}
\hline
\hline
Curve   &First\;\;appearance           &m      &p_1,\dots,p_6  &Group    &z-vector\\\hline
1-1     &66_{11}(twice)         &1      &3,0,1,0,0,0    &D_{3h} &36_{3,0}^{4}, 54_{3,0}\\
1-2     &126_{76}               &7      &0,3,4,3,0,0    &D_{3h} &18^{6}; 90_{9,0}^{3}\\
1-3     &106_{53}               &3      &1,2,3,0,0,0    &C_{s}  &18^{4}; 30_{1,0}, 42_{1,0}, 54_{1,0}, 60_{5,0}^{2}\\
1-4     &86_{39}                &4      &1,2,3,1,0,0    &C_1    &60_{6,0}^{2}, 138_{33,0}\\
1-5     &114_{36}(twice)        &5      &0,4,3,0,1,0    &C_{2v} &54_{3,0}, 72_{7,0}^{4}\\
1-6     &108_{48}(twice)        &7      &0,2,7,0,1,0    &C_{2v} &78_{9,0}^{3}, 90_{9,0}\\
1-7     &142_{82}               &5      &0,3,4,2,0,0    &C_1    &24^{3}; 60_{1,0}, 90_{8,0}^{2}, 114_{14,0}\\
1-8     &114_{61}               &11     &0,2,8,5,0,0    &C_{1}  &102_{13,0}^{2}, 138_{19,0}\\
2-1     &140_{68}               &2      &2,2,2,0,0,0    &C_2    &24^{2}; 72_{6,0}^{2}, 114_{13,0}^{2}\\
2-2     &90_{30}                &5      &0,6,0,3,0,0    &D_{3h} &18^{4}; 54_{3,0}, 72_{9,0}^{2}\\
2-3     &122_{21}               &5      &2,2,2,3,0,0    &C_{2}  &42_{1,0}, 78_{5,0}^{2}, 84_{9,0}^{2}\\
2-4     &126_{39}               &6      &0,4,4,2,0,0    &C_2    &36_{1,0}^{2}, 90_{10,0}^{2}, 126_{21,0}\\
2-5     &134_{130}              &6      &0,3,6,1,0,0    &C_{2v} &36^{2}; 42_{2,0}^{2}, 54_{3,0}, 96_{11,0}^{2}\\
2-6     &122_{104}              &7      &0,4,4,3,0,0    &C_{2v} &30^{4}; 54_{3,0}, 96_{11,0}^{2}\\
2-7     &124_{100}              &7      &0,3,6,2,0,0    &C_2    &90_{9,0}^{2}, 96_{11,0}^{2}\\
2-8     &128_{64}               &7      &0,2,6,3,0,0    &C_s    &54_{3,0}, 96_{11,0}^{2}, 138_{17,0}\\
%2-9    &126_{75}               &7      &0,3,6,2,0,0    &C_2    &30^{2}; 96_{11,0}^{2}, 126_{21,0}\\
2-9     &110_{24}               &8      &0,4,4,4,0,0    &C_2    &90_{12,0}^{2}, 150_{33,0}\\
2-10    &110_{74}               &8      &0,2,8,2,0,0    &C_2    &90_{12,0}^{2}, 150_{33,0}\\
%2-12   &134_{117}              &8      &0,2,8,2,0,0    &C_2    &102_{12,0}^{2}, 198_{45,0}\\
2-11    &134_{50}               &10     &0,4,4,6,0,0    &C_2    &114_{14,0}^{2}, 174_{29,0}\\
2-12    &134_{149}              &14     &0,2,8,8,0,0    &C_2    &126_{18,0}^{2}, 150_{21,0}\\
2-13    &126_{57}               &14     &0,0,12,6,0,0   &D_3    &126_{9,0}, 126_{18,0}^{2}\\
4-1     &72_{17}                &4      &0,6,4,0,0,0    &D_2    &24; 48_{4,0}, 72_{12,0}^{2}\\
6-1     &108_{18}               &6      &0,6,8,0,0,0    &D_3    &24^{3}, 36; 108_{18,0}^{2}\\
6-2     &138_{102}              &21     &0,1,20,6,2,0   &C_2    &78_{3,0}, 168_{33,0}^{2}\\\hline
6-3     &158_{150}              &24     &0,0,24,6,0,2   &D_3    &102_{9,0}, 186_{36,0}^{2}\\
8-1     &144_{151}              &8      &0,6,12,0,0,0   &D_2    &36;108_{12,0}, 144_{24,0}^{2}\\
\hline
\hline
\end{array}
\end{equation*}
}
\centering
\epsfxsize=120mm
\epsfbox{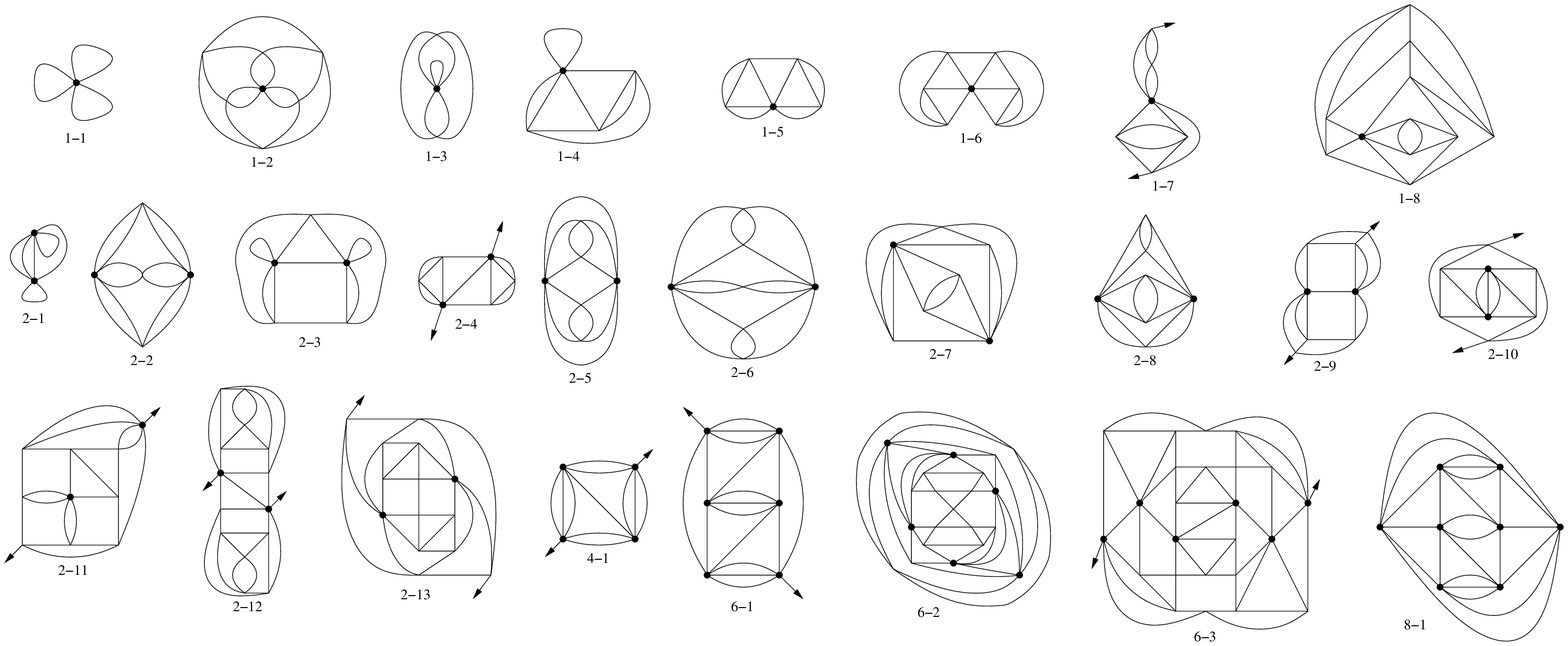}
\caption{Curves with triple self-intersections in railroads of graphs of type $4_n$, $n\leq 142$.}
\label{tab:Curve-with-triple-points}
\end{table}

\begin{remark}
Let $Z$ be a zigzag of signature $(\alpha_1,\alpha_2)$ bounding a
railroad $R$ with $m_2$ double and $m_3$ triple 
self-intersections. Then one has $m_2+3m_3=\alpha_1+\alpha_2$.
\end{remark}

\begin{proposition}\label{crude-upper-bound}
Let $G$ be a simple polyhedron with $p$-vector $p=(...,p_i,...)$ and let $G$ be tight. Then the number of zigzags of $G$ is at most 
$$\sum_{i\not= 6} ip_i /2\;\,.$$
\end{proposition}
\proof In fact, each zigzag has a non $6$-gonal face on each of its sides, since, otherwise, it would have a railroad on this side; so, the total number of incidences between zigzags and non $6$-gons is at least twice the number of zigzags. On the other hand, the number of those incidences is exactly $\sum_{i\not= 6} ip_i$. \qed

\begin{proposition}
Let $G$ be a graph of type $4_n$, having railroads $R_1$, \dots, $R_{p}$; let
$H$ be a plane graph formed by the curves representing those 
railroads. Then:

(i) Every $t$-gonal face of $H$ with $t=0,1$ contains exactly $3$-t $4$-gons of $G$; every $2$-gonal face of $H$ contains one or two $4$-gons of $G$.

(ii) If $q_t$ is the number of $t$-gonal faces of $H$, then one has the inequality
$$3q_0+2q_1+q_2\leq 6\;.$$ 

\end{proposition}
\proof Every $t$-gonal face $F$ of $H$ can be viewed as a regular patch with $t=t_{ob}+t_{ac}$ (obtuse and acute) intersections of arcs. Let $p'_4$ be the number of $4$-gons inside $F$. We will apply Theorem \ref{Local-Euler-Formula}.

(i) If $t=0$, then $2p'_4=6$ and we are done. If $t=1$, then
$$2p'_4=6-t_{ob}-2t_{ac}\geq 4\mbox{~~and~~}2p'_4=6-t_{ob}-2t_{ac}<6\;.$$
Since $2p'_4$ is even, we get that $2p'_4=4$, i.e., every $1$-gon contains exactly two $4$-gons. If $t=2$, then
$$2p'_4=6-t_{ob}-2t_{ac}\geq 2\mbox{~~and~~}2p'_4=6-t_{ob}-2t_{ac}<6\;.$$
This yields $p'_4=1$ or $p'_4=2$.

(ii) Any graph of type $4_n$ has six $4$-gons; so, the result follows. \qed

\bigskip

The Euler formula for the $p$-vector of a $3$-valent polyhedron, 
$12=\sum_{t} (6-t)p_t$, is a discrete analog of the Gauss-Bonnet 
formula, $2\pi(1-g)=\int_{S}K(x)dx$, for the Gaussian curvature 
$K$ of a surface $S$ of genus $g$. So, the $q$-gons
can be seen as {\em positively curved},
{\em flat}, or {\em negatively curved}, if $q<6$, $q=6$, or $q>6$, respectively. 
%See \cite{T} for deep development of those analogies.

Consider a two-faced graph $G$, which is a $3_n$, $4_n$ or $5_n$; it
has four $3$-gons, six $4$-gons, or twelve $5$-gons, respectively. Let
us call the {\em graph of curvatures} of $G$ the following graph (possibly,
with loops and multiple edges) having as vertex-set all non $6$-gonal
faces of $G$. Two vertices, say, non $6$-gonal faces $b$ and $c$ of
$G$ are called adjacent if there exist a pseudo-road
connecting faces $b$ and $c$. A {\em pseudo-road} is a sequence of
hexagons, say, $a_1$, \dots, $a_l$,
such that putting $a_0=b$ and $a_{l+1}=c$, we have that any $a_i$,
$1\leq i\leq l$, is adjacent to $a_{i-1}$ and $a_{i+1}$ on
opposite edges. Clearly, the graph of curvatures of a
graph $w_n$, $w\in \{3,4,5\}$, is regular of degree $w$.

Given an hexagon in a graph of type $3_n$, $4_n$, $5_n$, any pair
of opposite edges belongs to a railroad or a pseudo-road; so,
we exactly have a triple covering of the set of all hexagons
by the set of all railroads and pseudo-roads. Every non $6$-gonal
face, i.e., $t$-gonal face with $t<6$, has $t$ adjacencies with
the system of pseudo-roads.

\begin{figure}
\centering
\begin{minipage}[b]{7cm}%
\leavevmode
\begin{center}
\epsfig{figure=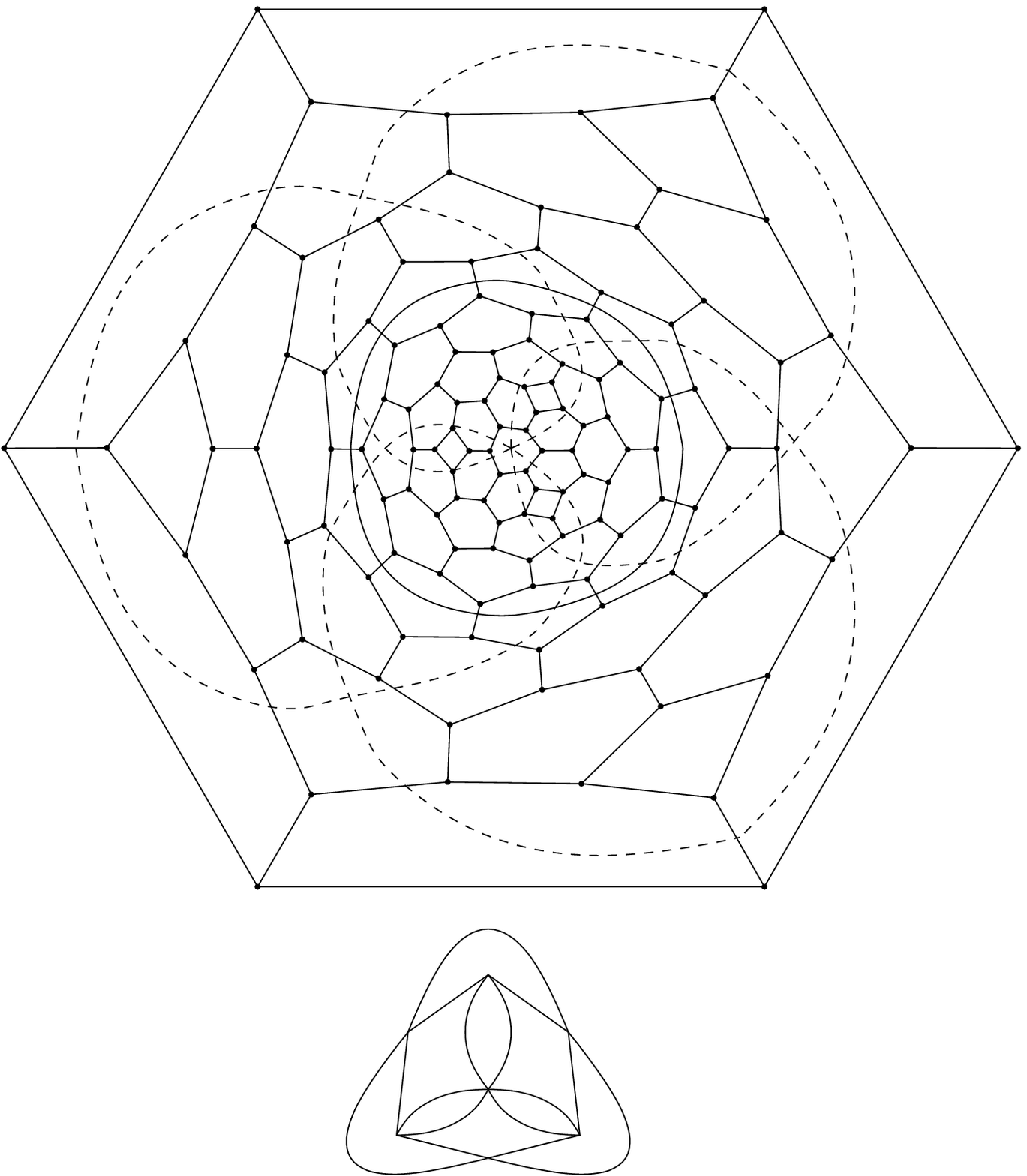,height=4.5cm}
\end{center}
\begin{center}
One of two self-intersecting railroads and the equatorial simple railroad
\end{center}
\end{minipage}
\begin{minipage}[b]{6.5cm}%
\leavevmode
\begin{center}
\epsfig{figure=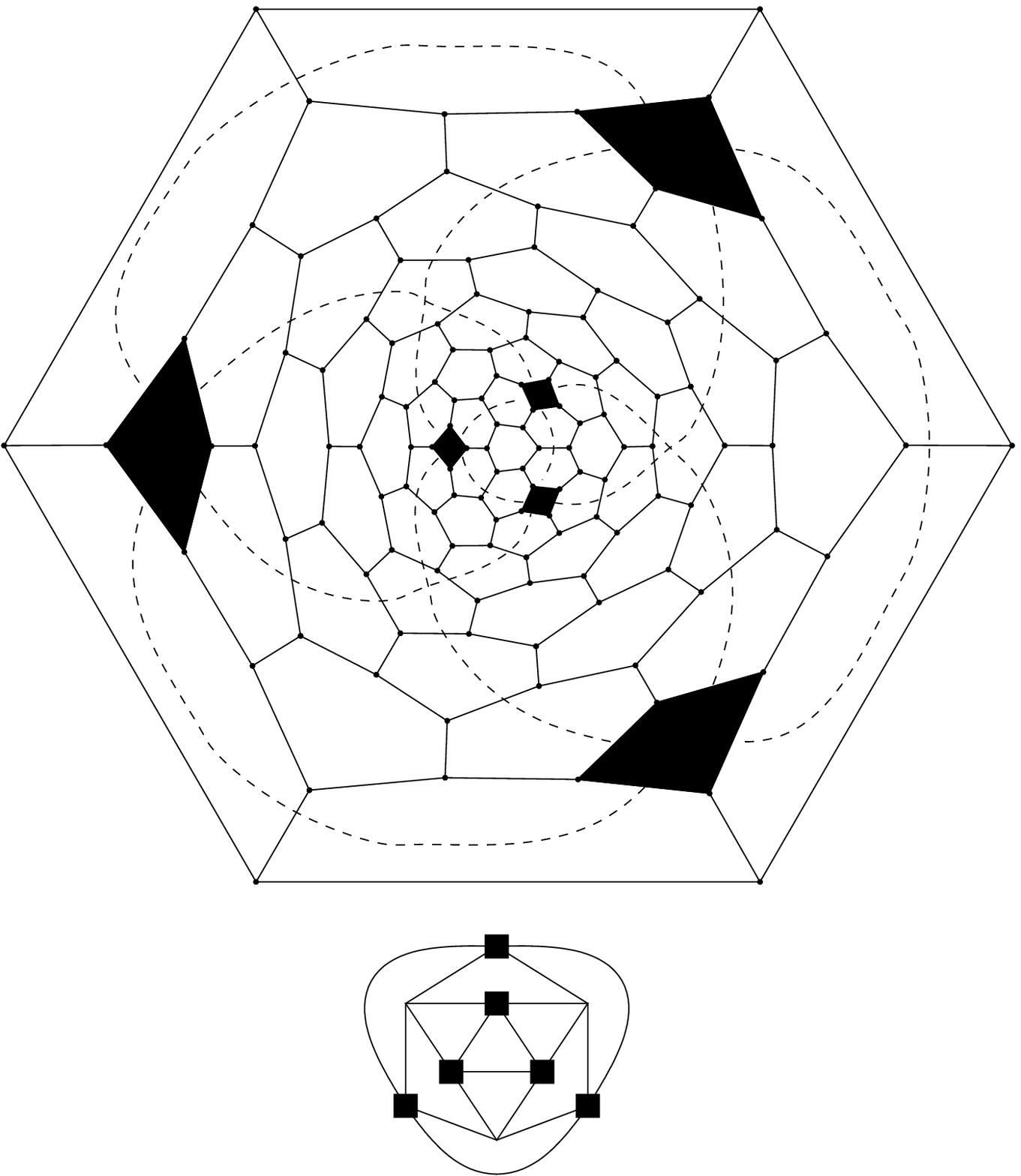,height=5cm}
\end{center}
\begin{center}
All twelve pseudo-roads
\end{center}
\end{minipage}
\caption{Road-decomposition of a graph $4_{126}(D_{3h})$.}
\label{fig:RailRoadSystem}
\end{figure}

For the special case, when our graph is of type $4_n$, any pseudo-road 
arriving on a $4$-gonal face can be extended on the opposite edge. So, the 
set of such extended pseudo-roads, together with the set of railroads,
can be identified with the set of 
central circuits of the dual of our graph. Hence, extended pseudo-roads 
can be seen (in the same way as it was done for railroads) as
projections of a Jordan curve on the plane with double 
and, eventually, triple points of self-intersection. Those notions are
illustrated in Figure \ref{fig:RailRoadSystem}. This
graph $4_{126}(D_{3h})$ has the following 
{\em road decomposition}: five concentric simple railroads 
(we indicate only the central, {\em equatorial} one), two 
self-intersecting railroads 
(they differ only by their opposite position on the sphere; we present 
only one of them) and all $12$ pseudo-roads.

%\quad
%\subfigure[Graphs]{\epsfxsize=25mm\epsfbox{RailRoadSystemDraws.eps}}}

\section{Two-faced polyhedra}
Remind that {\em two-faced} polyhedra have only $a$- and $b$-gonal 
faces with $3\leq a < b\leq 6$ (see \cite{DG2}). Denote 
by $p_a>0$ and $p_b\geq 0$ the number of $a$-gonal and $b$-gonal 
faces.

Any $3$-valent $n$-vertex two-faced polyhedron has
$3n/2=(ap_{a}+bp_{b})/2$ edges and satisfies the Euler relation
$n-3 \frac{n}{2}+(p_{a}+p_{b})=2$, i.e., $n=2(p_{a}+p_{b}-2)$ and
$p_{a}(6-a)+p_{b}(6-b)=12$. 
So, $b<6$ is possible only for $6$ simple polyhedra 
with $(a,b)=(3,4)$($Prism_3$), $(3,5)$(D\"urer's Octahedron), $(4,5)$(four dual deltahedra). The other cases, namely, 
$(a,b)=(3,6),(4,6),(5,6)$ are denoted by $3_n$, $4_n$, $5_n$, 
where $n$ is the number of vertices. The polyhedra $5_n$ are the well known 
{\em fullerenes} (of Organic Chemistry).
See Table \ref{tab:listOfCases} for all the possibilities.
The criterion for the existence of the polyhedra of type $3_n$, $4_n$
and $5_n$ is due to \cite{GrMo}.

\begin{table}
\[\begin{array}{||c||c|c|c|c||}
\hline
\hline
(a,b) & Polyhedra & \mbox{Exist~if~and~only~if} & p_a & n\\ \hline\hline
(5,6) & 5_{n} \mbox{~(fullerenes)} & p_6\in N-\{1\}& p_5=12 &n=20+2p_6\\ \hline
(4,6) & 4_{n} & p_{6} \in N-\{1\} & p_4=6 &n=8+2p_6\\ \hline
(3,6) & 3_{n} & p_{6}/2 \in N-\{1\} & p_3=4 &4+2p_6\\  \hline\hline
(4,5) & \mbox{4 dual deltahedra} & p_5=2,3,4,5 &p_4=5,4,3,2&n=10,12,14,16\\ \hline
(3,5) & \mbox{D\"urer's Octahedron} & p_5=6 &p_3=2 &n=12\\ \hline
(3,4) & Prism_{3} & p_4=3 & p_3=2 & n=6\\ 
\hline
\hline
\end{array}\]
\caption{All simple polyhedra with only $a$-gonal and $b$-gonal faces, $3 \le a< b \le 6$.}
\label{tab:listOfCases}
\end{table}

\begin{figure}
\begin{center}
\epsfxsize=150mm
\epsfbox{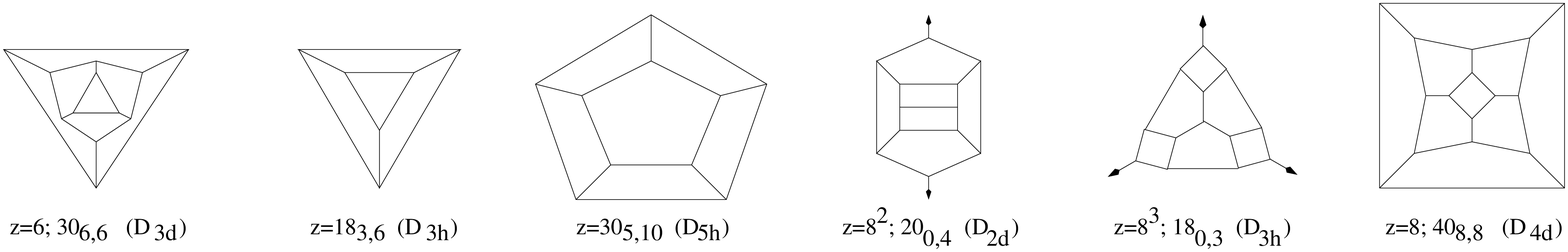}
\end{center}
\caption{D\"urer's Octahedron, $Prism_3$ and four dual deltahedra.}
\label{fig:Deltahedra}
\end{figure}

For the types $3_n$, $4_n$ and $5_n$ in Table \ref{tab:listOfCases}
the case $p_{6}=0$ yields the Tetrahedron, the Cube and the Dodecahedron.
Those three polyhedra, together with $4$ dual 
deltahedra and $Prism_3$, are the duals of the eight {\em convex deltahedra}.
{\em D\"urer's Octahedron} is the Cube truncated at two opposite vertices.
The four dual deltahedra from the fourth line of Table \ref{tab:listOfCases}
are: $Prism_5$ (dual $5$-bipyramid), dual {\em bisdisphenoid},
dual $3$-{\em augmented} $Prism_{3}$, dual $2$-{\em capped} $APrism_4$.
Those four dual deltahedra, preceded by the D\"urer's Octahedron
and $Prism_3$, are given in Figure \ref{fig:Deltahedra}
with their $z$-vectors and groups.

\begin{lemma}
Any connected $3$-valent plane graph, having only $q$-gonal faces with $3\leq q\leq 6$, is $2$-connected.
\end{lemma}
\proof Let $G$ be one such graph and assume that there is one vertex $v$, such that $G-\{v\}$ is disconnected. $G-\{v\}$ has two or three components; clearly, it cannot have three components, since $q\leq 6$.
Denote by $C_1$ and $C_2$ the two components and by $e=\{v,v'\}$ the edge
linking $C_1$ to $C_2$. Then two edges from $v$ will connect to another
vertex $w$ and two edges from $v'$ will connect to another vertex $w'$,
since we assume that the faces are incident to at most $6$ edges.
See below the corresponding drawing.

\begin{center}
\epsfxsize=40mm
\epsfbox{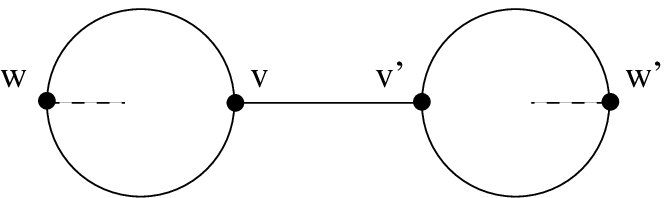}
\end{center}

But the vertex $w$ will disconnect the graph and so, iterating the
construction, we obtain an infinite sequence $v_1$, \dots, $v_n$ of
vertices that disconnect $G$. This contradicts to initial assumption and proves that $G$ is $2$-connected. \qed

\begin{figure}
\begin{center}
\epsfxsize=60mm
\epsfbox{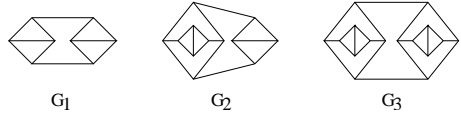}
\end{center}
\caption{First graphs $(G_n)_{n\geq 1}$.}
\label{SequenceOfGraphs}
\end{figure}

\bigskip
Denote by $(G_n)_{n\geq 1}$ the $3$-valent plane 
$4(n+1)$-vertex graph, whose faces are (organized in pairs of
adjacent ones) triangles and hexagons.
The graph $G_n$ is $2$-connected but not $3$-connected, its $z$-vector 
is $4^{n+1}, (4n+4)^2$; its symmetry group is $D_{2d}$ or $D_{2h}$,
if $n$ is even or odd, respectively. The first occurrences are depicted 
in Figure \ref{SequenceOfGraphs}.

\begin{figure}
\centering
\begin{minipage}[b]{4cm}%
\centering
\epsfig{figure=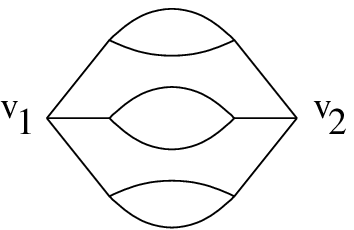,height=2.5cm}\par
a) 
\end{minipage}
\begin{minipage}[b]{4cm}%
\centering
\epsfig{figure=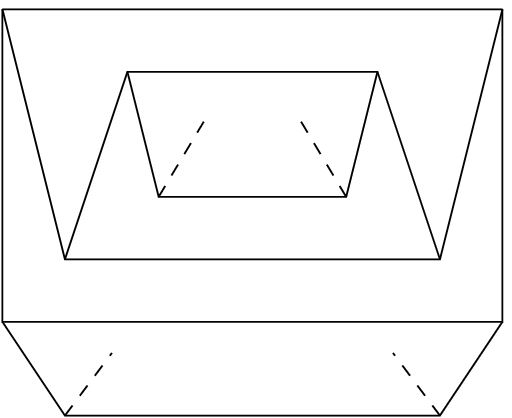,height=2.5cm}\par
b) 
\end{minipage}
\begin{minipage}[b]{4cm}%
\centering
\epsfig{figure=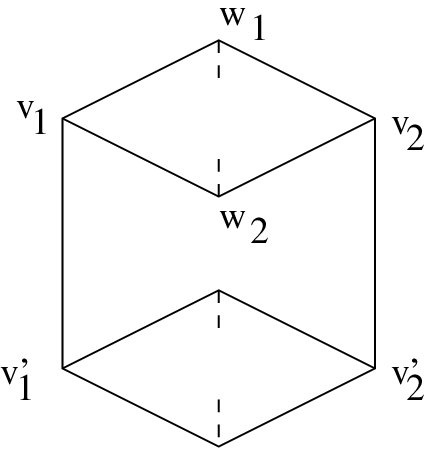,height=2.5cm}\par
c) 
\end{minipage}
\caption{Some $2$-connected graphs}
\label{StepProof}
\end{figure}

\begin{proposition}\label{3-connectedness}
Any $3$-valent plane graph having only $q$-gonal faces with $3\leq q\leq 6$, is
$3$-connected, except the graphs of the family $(G_n)_{n\geq 1}$.
\end{proposition}
\proof Let $G$ be a $3$-valent plane graph with $k$-gonal faces, $3\leq k\leq 6$ and assume that it is not $3$-connected. Then there are two
vertices, say, $v_1$ and $v_2$, such that $G-\{v_1, v_2\}$ is 
disconnected. If $G-\{v_1,v_2\}$ has three components, then
since $q\leq 6$, it is of the form depicted in Figure \ref{StepProof}.a), which is impossible by the condition $q\geq 3$; so, $G-\{v_1,v_2\}$ has two components, say, $C_1$ and $C_2$. 

There are two edges, say, $e_1=\{v_1, v'_1\}$ and $e_2=\{v_2, v'_2\}$, that connect $C_1$ to $C_2$. Since $G-\{v_1, v_2\}$ is disconnected, $e_1$ and $e_2$ are both incident to faces $F$ and $F'$. Then, by our assumption on size of faces, we get the following possibilities:

\begin{center}
\epsfxsize=145mm
\epsfbox{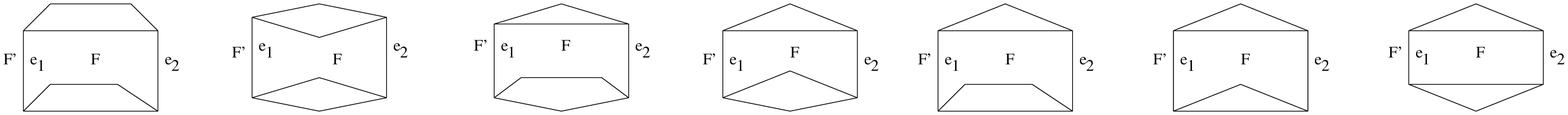}
\end{center}

Only first two cases are possible, since, otherwise, we would have a non
$2$-connected plane graph. See in Figure \ref{StepProof}.b),c) the possible
continuation of those graphs.

In case b) one obtain an infinite structure; so, we get again a
contradiction. 

Let us consider now the case c); the two points, $w_1$ and $w_2$, can,
either be connected by an edge and we are done, or be non-connected, in 
which case $w_1$ and $w_2$ disconnect the
graph. In the latter case, one can iterate the construction. Since
the graph is finite, the construction eventually finish and we get a 
graph $G_{n_0}$, $n_0\geq 1$. \qed

\bigskip

All graphs of type $3_n$, $4_n$, $5_n$ with maximal pair of symmetry
($T$ or $T_d$, $O$ or $O_h$, $I$ or $I_h$, respectively) are known;
in fact, the Goldberg-Coxeter construction, denoted $GC_{k,l}$
(\cite{Cox71}, \cite{Gold}, \cite{DD03}),
implies that any graph of this type $3_n$, $4_n$, $5_n$ and those
symmetry are of the form
$GC_{k,l}(G_0)$ with $0\leq l\leq k$ and $G_0$ being Tetrahedron, Cube
or Dodecahedron, respectively.
Those graphs have $n_0(k^2+kl+l^2)$ vertices with $n_0=4$, $8$ or
$20$, respectively. Those graphs are tight if and only if $gcd(k,l)=1$;
they are of symmetry $T_{d}$, $O_h$ or $I_h$ if and only if $l=0$ or $k=l$.

% with
%$k\geq l\geq 0$ and its symmetry is $T_d$ if and only if $l=0$ or
%$k=l$. It is easy to see that in the first case (i.e., if $n=4k^2$)
%$z(G)=(4k)^{3k}$ and $Int=(2^{2k},0^{k-1})$, while in the second case
%(i.e., $n=12k^2$) we have $z(G)=(12k)^{3k}$ and
%$Int=(6^{2k},0^{k-1})$. The case of tight $3_n$ of symmetry $T$ or
%$T_d$ corresponds to $gcd(k,l)=1$, and the only tight $3_n$ of symmetry
%$T_d$ are the Tetrahedron and the Truncated Tetrahedron.

The polyhedra dual to the $3$-valent polyhedra without $q$-gonal faces,
$q>6$, are studied in \cite{T}; they are
called there {\em non-negatively curved triangulations} (actually, 
the reference \cite{sah94} is an application of a preliminary 
version of \cite{T}). 
Theorem 3.4 in \cite{sah94} implies that $N_3(n)=O(n)$, $N_4(n)=O(n^3)$ 
and $N_5(n)=O(n^9)$, where $N_w(n)$ denotes the number of polyhedra
$w_n$ for $w=3, 4, 5$, respectively.

%The main Theorem 0.1 in \cite{T} describes 
%them as the elements of $L_+/G$, where $L$ is a lattice in 
%complex Lorenz space $C^{(1,9)}$, $G$ is a group of automorphisms 
%and $L_+$ is the set of lattice points 
%of positive square-norm. The square of the norm of a lattice point 
%is the number of triangles in the triangulation and the number of 
%such triangulations with up to $2n$ triangles, is $O(n^{10})$.

\section{Polyhedra $3_n$}

We present in Table \ref{tab:listOf3nGraphs} the numbers $N_3(n)$ of
$3_n$ and the numbers $N^t_3(n)$ of tight $3_n$ for $n\leq 512$.

\begin{table}
\scriptsize
\begin{equation*}
\begin{array}{||c|c|c||c|c|c||c|c|c||c|c|c||c|c|c||c|c|c||c|c|c||}
\hline
\hline
n&N_{3}&N^t_3&n&N_{3}&N^t_3&n&N_{3}&N^t_3&n&N_{3}&N^t_3&n&N_{3}&N^t_3&n&N_{3}&N^t_3&n&N_{3}&N^t_3\\\hline
12&1&1&84&7&2&156&11&3&228&15&4&300&23&3&372&23&6&444&27&7\\
16&2&0&88&6&0&160&19&0&232&15&0&304&26&0&376&24&0&448&48&0\\
20&1&1&92&4&4&164&7&7&236&10&10&308&17&8&380&21&9&452&19&19\\
24&2&0&96&14&0&168&17&0&240&33&0&312&29&0&384&50&0&456&41&0\\
28&2&2&100&6&3&172&8&8&244&11&11&316&14&14&388&17&17&460&25&11\\
32&4&0&104&7&0&176&16&0&248&16&0&320&37&0&392&29&0&464&37&0\\
36&3&1&108&8&2&180&15&2&252&20&3&324&22&5&396&28&5&468&33&6\\
40&3&0&112&12&0&184&12&0&256&26&0&328&21&0&400&40&0&472&30&0\\
44&2&2&116&5&5&188&8&8&260&15&6&332&14&14&404&17&17&476&25&13\\
48&7&0&120&13&0&192&27&0&264&25&0&336&43&0&408&37&0&480&69&0\\
52&3&3&124&6&6&196&11&7&268&12&12&340&19&8&412&18&18&484&23&17\\
56&4&0&128&14&0&200&16&0&272&23&0&344&22&0&416&39&0&488&31&0\\
60&5&1&132&9&2&204&13&3&276&17&4&348&21&5&420&35&3&492&29&7\\
64&8&0&136&9&0&208&19&0&280&25&0&352&34&0&424&27&0&496&40&0\\
68&3&3&140&9&3&212&9&9&284&12&12&356&15&15&428&18&18&500&27&13\\
72&7&0&144&19&0&216&21&0&288&39&0&360&41&0&432&52&0&504&54&0\\
76&4&4&148&7&7&220&13&5&292&13&13&364&21&11&436&19&19&508&22&22\\
80&9&0&152&10&0&224&24&0&296&19&0&368&30&0&440&37&0&512&48&0\\
\hline
\hline
\end{array}
\end{equation*}
\caption{Numbers $N_3(n)$ of graphs $3_n$ and $N_3^t(n)$ of tight $3_n$ for $n\leq 512$.}
\label{tab:listOf3nGraphs}
\end{table}

\begin{figure}
\centering
\epsfxsize=60mm
\epsfbox{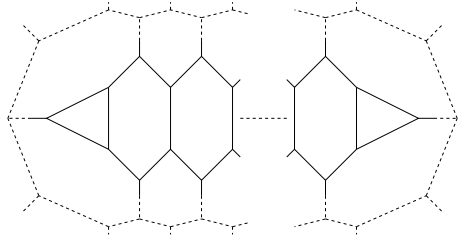}
\caption{A pseudo-road in a graph of type $3_n$ and a circular railroad.}
\label{fig:GrunbaumConstruction}
\end{figure}

\begin{figure}
\centering
\begin{minipage}[b]{3.2cm}%
\centering
\epsfig{figure=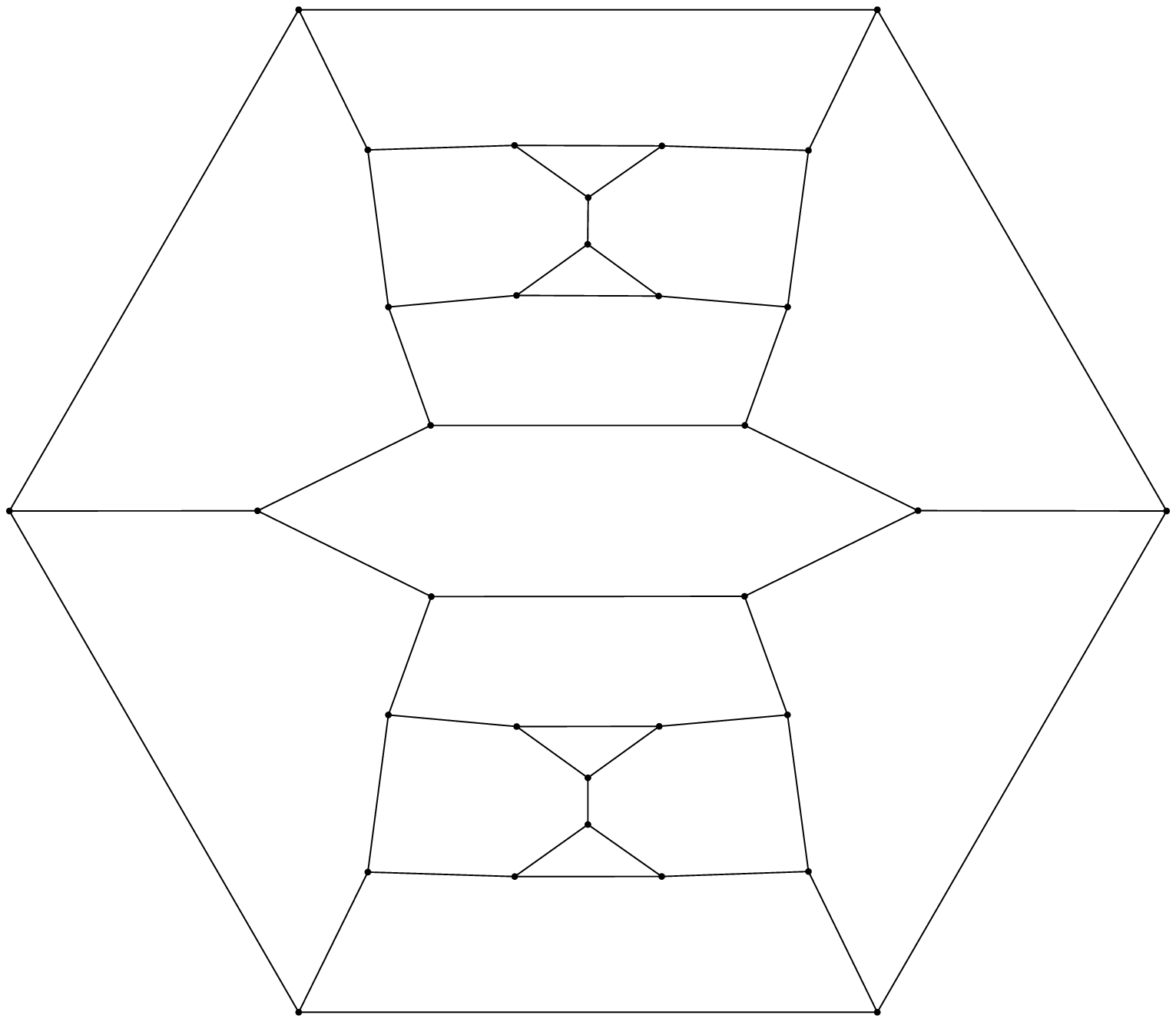,width=3cm}
a) $32$ vertices, case (iii), $D_{2h}$
\end{minipage}
\begin{minipage}[b]{3.2cm}%
\centering
\epsfig{figure=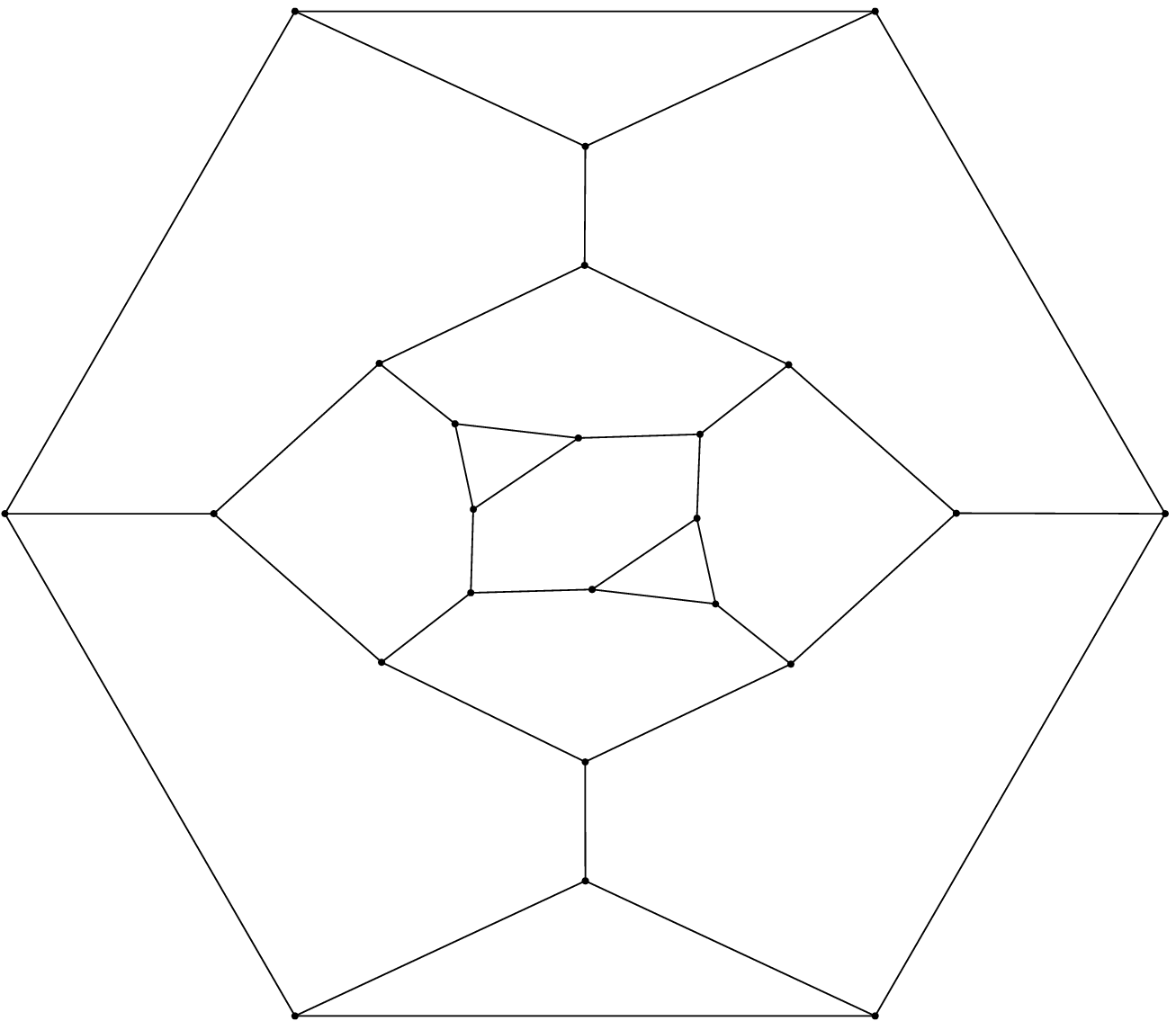,width=3cm}
b) $24$ vertices, case (ii), $D_{2}$
\end{minipage}
\begin{minipage}[b]{3.2cm}%
\centering
\epsfig{figure=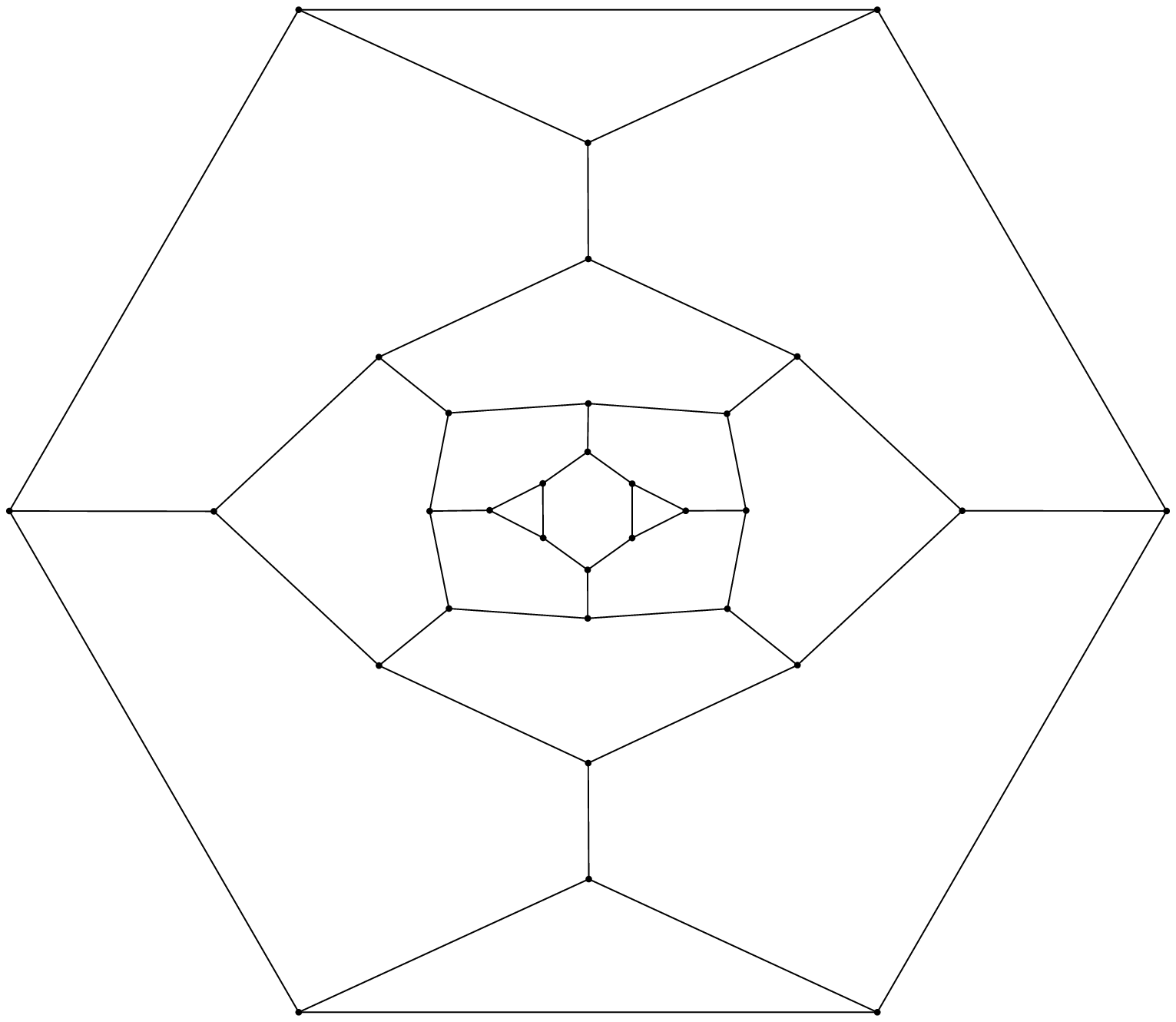,width=3cm}
c) $32$ vertices, case (ii), $D_{2d}$
\end{minipage}
\begin{minipage}[b]{3.2cm}%
\centering
\epsfig{figure=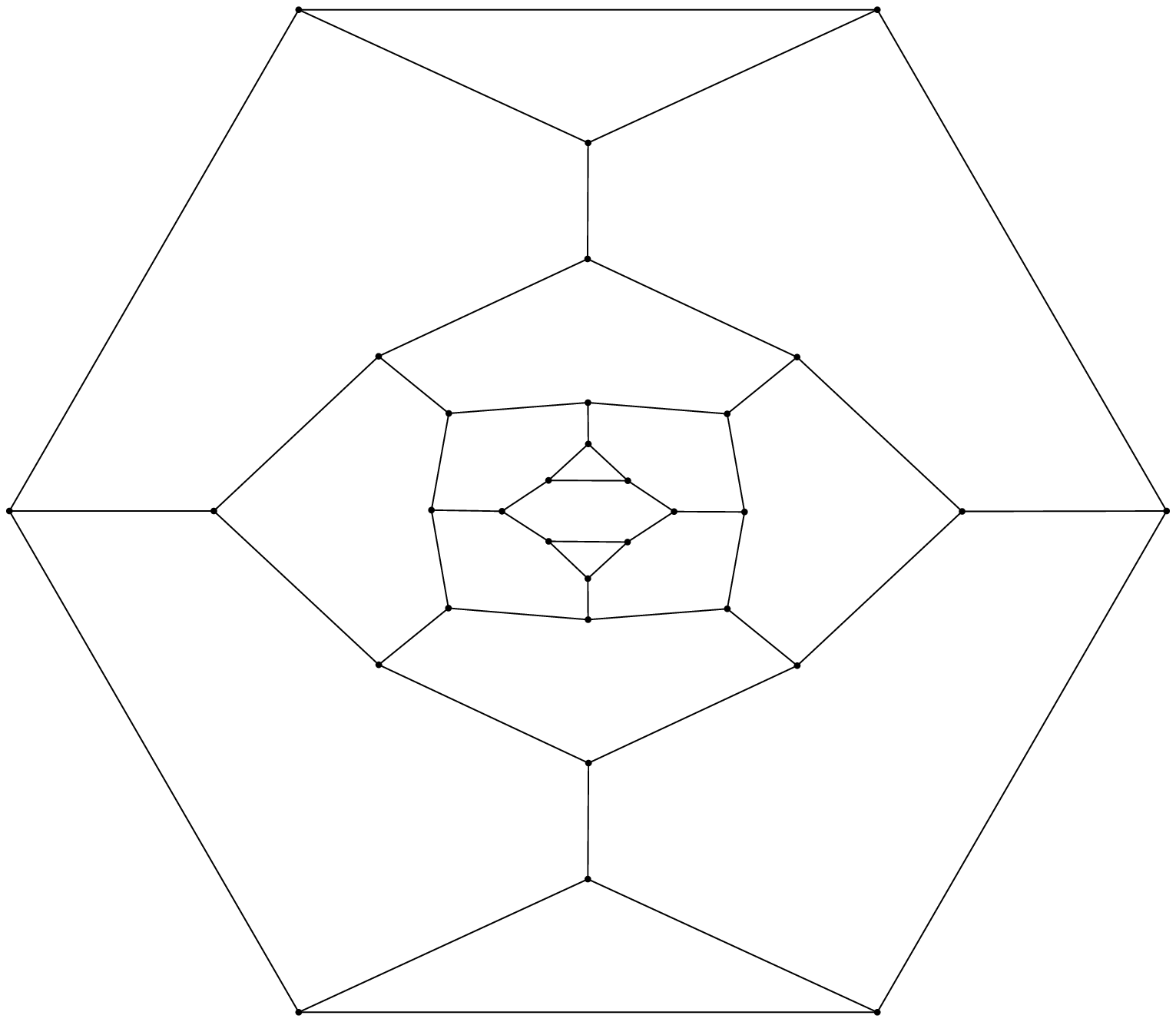,width=3cm}
d) $32$ vertices, case (ii), $D_{2h}$
\end{minipage}
\caption{Graphs with hexagons adjacent to two triangles.}
\label{fig:Case32}
\end{figure}

In \cite{GrMo} it was shown that each zigzag in a graph of type $3_n$ is simple and a general construction of all $3_n$ was given; we sketch this construction with slightly different notation.
If a zigzag contains two edges, say, $e_1$ and $e_2$, of a triangle, then
its third edge, say, $e_3$ defines a pseudo-road $PR$ that finish on another
triangle.
Denote by $s$ one fourth of the length of the zigzag; then $s-1$ is
the number of hexagons of $PR$. Consider the sequence of concentric 
simple railroads,
which, possibly, goes around the patch depicted in Figure
\ref{fig:GrunbaumConstruction}; let $m$ denote the number of zigzags and
$m-1$ be the number of corresponding
railroads. The same patch (with the same number $s-1$ of hexagons
between the other two triangles and with the same number $m-1$ of concentric
simple railroads, going around it) occurs on the other side of the
sphere. As in \cite{GrMo}, we get, by direct computation, the equality 
$p_6=2(sm-1)$ and $n=4sm$, where $p_6$ is the number of hexagons. See 
Figure \ref{fig:Case32}.c),d) for two examples with $s=2$ and $m=4$.

Take a triangle $T$ and an edge $e$ of this triangle; $e$ belongs to a sequence of adjacent $6$-gons, which is concluded by another triangle $T'$. If one considers the three edges $e_1$, $e_2$, $e_3$, then we obtain three triangles $T_1$, $T_2$, $T_3$, respectively, that may be distinct or not.

\begin{theorem}\label{Possible-forms-for-3nS}
Any graph $3_n$ has exactly one of the following forms:

(i) The Tetrahedron or the Truncated Tetrahedron
(the only case, when an hexagon is adjacent to more than two triangles).

(ii) There are two hexagons, every of which is adjacent to two
triangles on opposite edges; there are one or two such graphs of type $3_n$
depending on $n\equiv 8\pmod {16}$ or $n\equiv 0 \pmod {16}$,
respectively. Their symmetry groups are, respectively, $D_2$ if
$n\equiv 8\pmod {16}$ and $D_{2h}$, $D_{2d}$ if $n\equiv 0\pmod {16}$.

(iii) There are four hexagons (in two adjacent pairs), every of which
is adjacent to two triangles on non-opposite and non-adjacent edges;
there is exactly one such graph $3_n$ for every $n\geq 16$, $n\equiv 0
\pmod 4$.\footnote{In the exceptional case $n=16$, there are only two
polyhedra $3_{16}$: the Cube truncated at four non-adjacent vertices and
the Cube truncated at four vertices of two opposite edges. The
second polyhedron is the unique and only example, which is of type
(ii) and (iii).}
The symmetry group is $D_{2h}$ or $D_{2d}$ if $\frac{n}{4}$ is even or odd,
respectively.

(iv) Each hexagon is adjacent to at most one triangle; there
is an one-to-one correspondence between such graphs $3_n$, having
isolated triangles, and IPT fullerenes $5_n$ (i.e., those having their
$12$ pentagons organized in $4$ triples of mutually adjacent ones).
\end{theorem}
\proof For given graph $G$ of type $3_n$, denote by $t(G)$ the maximal
number of triangles, which are adjacent to a hexagon.
The only graph of type $3_n$ without hexagons is the Tetrahedron.
The case $t(G)=4$ corresponds to the $2$-connected but not $3$-connected
$6$-vertex plane graph $G_1$ (see Figure \ref{SequenceOfGraphs}).
The only graph of type $3_n$ with $t(G)=3$ is the Truncated Tetrahedron.
If $t(G)=2$, then two cases occur: either two triangles are adjacent to
an hexagon on opposite edges or they are adjacent to an hexagon on
non-opposite and non-adjacent edges.

The first case corresponds to the case $s=2$ of the Gr\"{u}nbaum-Motzkin
construction. The divisibility of $n$ by $8$ follows from the formula $n=4sm$.
At the last step of the Gr\"{u}nbaum-Motzkin construction, namely,
after adding railroads, we add a pair of triangles in one of two possible ways.
If $n\equiv 8\pmod {16}$, then both ways give isomorphic graphs, both of
symmetry $D_2$, for $n\geq 8$ (see Figure \ref{fig:Case32}.b), for this
possibility). If $n\equiv 0\pmod {16}$, then we get two non-isomorphic graphs,
one of symmetry $D_{2h}$, the other of symmetry $D_{2d}$ (see Figure
\ref{fig:Case32}.c),d) for those two possibilities).

Let our graph be of type (iii), i.e., two hexagons $H_1$, $H_2$ are adjacent to two triangles $T_1$, $T_2$ on non-opposite and non-adjacent edges, as in the picture below.

\begin{center}
\epsfxsize=60mm
\epsfbox{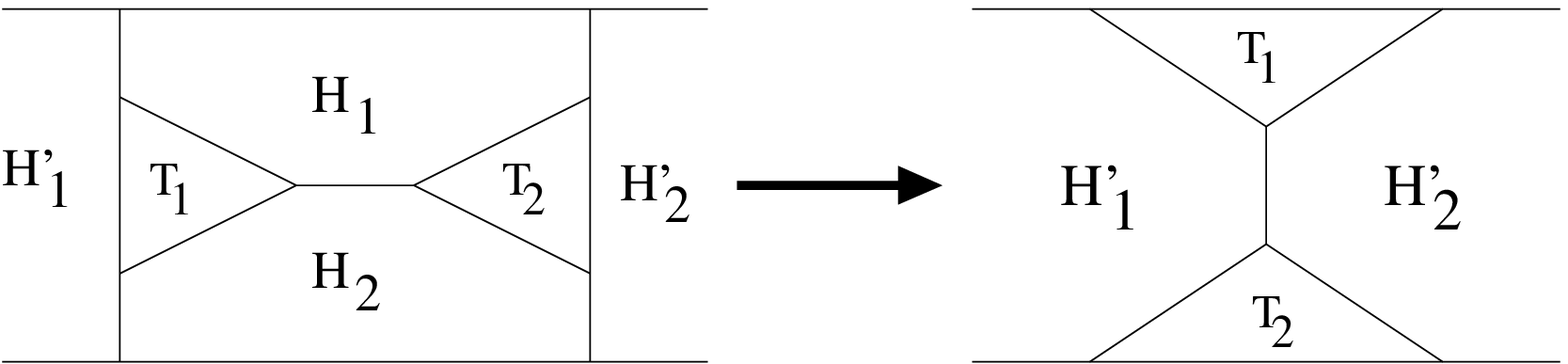}
\end{center}

These two hexagons $H_1$, $H_2$ are adjacent, respectively, to two hexagons $H'_1$, $H'_2$. By replacing the patch $H_1$, $H_2$, $T_1$, $T_2$ by two triangles, one gets a graph $3_{n-4}$, which is also of type (iii). By induction, one gets, for any $n\equiv 0\pmod 4$ with $n\geq 16$, an unique graph.
A graph of type (iii) has two pairs of hexagons, each pair being adjacent to two triangles; if one does the operation, depicted above, to the two pairs, then the symmetry group remains the same. So, one has a periodicity of order $8$ for the symmetry groups and the result on groups follows.

If the graph has $t(G)=1$, then the operation, depicted in picture below, on all triangles realize a one-to-one correspondence. \qed

\begin{center}
\epsfxsize=60mm
\epsfbox{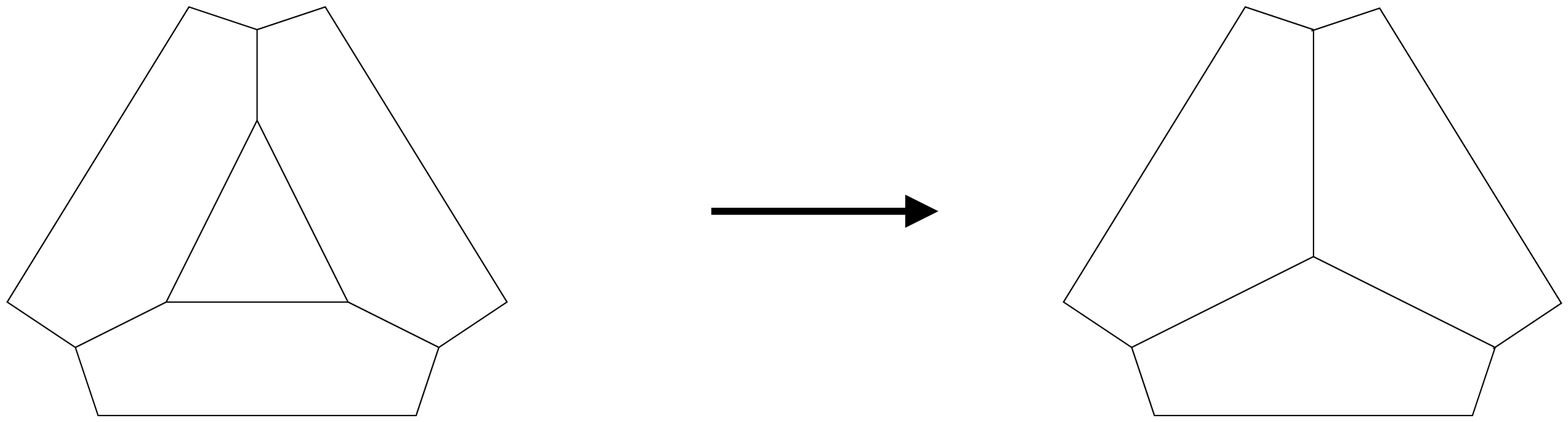}
\end{center}

See Figure \ref{fig:Case32}.a),c),d) for all graphs $3_{32}$ having hexagons adjacent to two triangles.

%NEED TO ELIMINATE (a), (b) (c) in subfigure environment

\begin{proposition}\label{PossibleGraphCurvature}
The graph of curvatures of a $3_n$ is one of the three following $4$-vertex graphs:
\begin{center}
\epsfxsize=80mm
\epsfbox{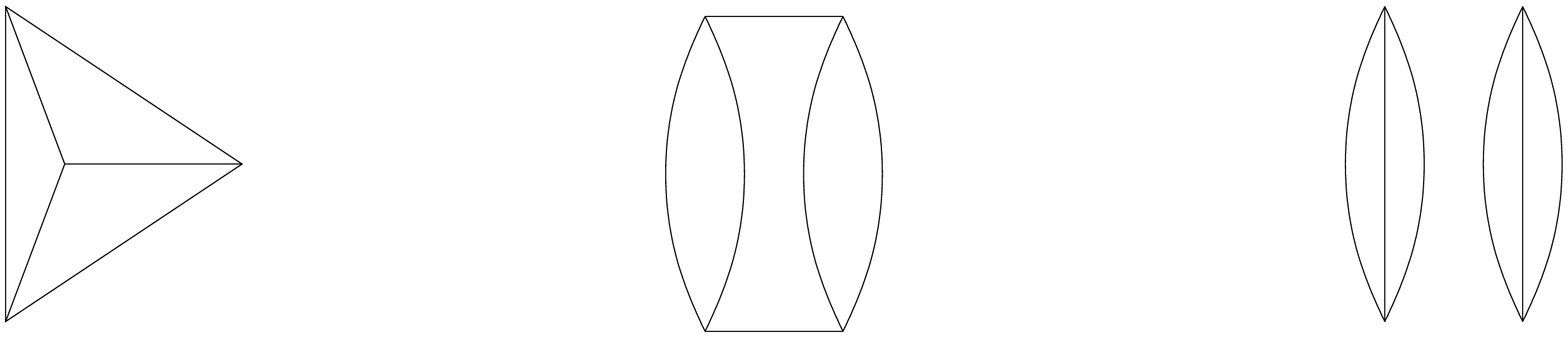}
\end{center}

\end{proposition}
\proof Take a triangle, say, $T_1$ and an edge $e_1$ of $T_1$. The
pseudo-road, which is defined by $T_1$ and $e_1$, establishes an edge
between $T_1$ and, say, $T_2$; this pseudo-road is bounded by one
zigzag $Z_1$. This zigzag belongs to a sequence of $m$ concentric zigzags
$Z_1, Z_2, \dots, Z_m$. The zigzag $Z_m$ defines a pseudo-road
between the remaining triangles $T_3$ and $T_4$, i.e., an
edge connecting $T_3$ to
$T_4$ in the graph of curvatures. Since any vertex of the
graph of curvatures has degree $3$, we are done. \qed

\begin{conjecture}
The graph of curvatures of any tight graph of type $3_n$
is as in the first case of Theorem \ref{PossibleGraphCurvature}.
\end{conjecture}
%This conjecture was checked for $n\leq 500$.

Fix a graph $G$ of type $3_n$.
Denote by $T_i$, $1\leq i\leq 4$, the triangles in $G$, by
$(e_i)_{1\leq i\leq 3}$ the edges of $T_1$, by $PR_i$ the pseudo-road
defined by $T_1$, $e_i$ and by $(4s_i)_{1\leq i\leq 3}$ the lengths of
corresponding zigzags. Without loss of generality, we can suppose that
$s_1\leq s_2\leq s_3$. Denote by $m_1-1$, $m_2-1$, $m_3-1$,
respectively, the corresponding number of concentric railroads around
each of three patches.
Denote by $Z_{i,k}$, where $1\leq i\leq 3$ and $1\leq k\leq m_i$, the zigzags of $G$.

\begin{theorem}\label{Theorem-interest-ZigZag-3n}
For a graph $G$ of type $3_n$, the following properties hold:

(i) Its $z$-vector is $(4s_1)^{m_1}, (4s_2)^{m_2},(4s_3)^{m_3}$; the number of railroads is $m_1+m_2+m_3-3$.

(ii) $G$ has at least three zigzags with equality if and only if it is tight.

(iii) If $G$ is tight, then $z(G)=n^3$ (so, each zigzag is a Hamiltonian circuit).

(iv) $G$ is $z$-balanced and $|Z_{i,k}\cap Z_{j,l}|=\left\lbrace\begin{array}{lcl}
0               &\mbox{~if~}    &i=j,\\
\frac{n}{2m_im_j}     &\mbox{~if~}    &i\not= j.
\end{array}\right.$.

\end{theorem}
\proof If $Z$ is a fixed zigzag of $G$, then it belongs to a sequence of $m_i$ concentric zigzags of a patch defined by two triangles and a pseudo-road. (i) follow immediately; (ii) is a corollary of (i), since $G$ is tight if and only if all $m_i$ are equal to $1$.

Let $G$ be a tight $3_n$; the formula $n=4sm$ becomes $n=4s_1=4s_2=4s_3$. So, $z=n^3$. 

Clearly, any two zigzags $Z_{i,k}$ and $Z_{j,l}$ are
disjoint if $i=j$. Moreover, the size of the intersection $Z_{i,k}\cap
Z_{j,l}$ depends only on $i$ and $j$. Denote by $\beta_{ij}$ the size of the pairwise intersections $|Z_{i,1}\cap Z_{j,1}|$. We obtain the linear system 
$$\left\lbrace\begin{array}{rcl}
4s_1&=&            \beta_{13} m_3+\beta_{12}m_2\\
4s_2&=&\beta_{23}m_3+             \beta_{12}m_1\\
4s_3&=&\beta_{23}m_2+\beta_{13}m_1
\end{array}\right. ,$$
which has the unique solution $\beta_{ij}=\frac{n}{2m_im_j}$. If one writes this intersection size as $8s_is_j/n$, then it is easy to see that $G$ is $z$-balanced. \qed

\begin{theorem}
All graphs of type $3_n$ are tight if and only if $\frac{n}{4}$ is prime.
%If $p\leq q$ are both prime then $N_{3}(4p)\leq N_3(4q)$
\end{theorem}
\proof Assume that $\frac{n}{4}$ is prime; let $G$ be a graph $3_n$ with
parameters $s_i$, $m_i$ defined before Theorem \ref{Theorem-interest-ZigZag-3n}. If $\frac{n}{4}=s_im_i$, then either $s_i=1$, or $m_i=1$. The first case corresponds to a $2$-connected but not $3$-connected plane graph; so, $m_1=m_2=m_3=1$ and $G$ has no railroad.
The second case $m_i=1$ means also the absence of a railroad.

Assume that $\frac{n}{4}$ is not prime; if $n=4pq$ for $q>1$, then, using
the Gr\"{u}nbaum-Motzkin construction, we can construct a graph of type
$3_n$ with a
system of $q>1$ concentric zigzags of length $4p$ and so, with
at least $q-1\geq 1$ railroads. \qed

\begin{remark}
In Table \ref{tab:listOf3nGraphs}, the number of graphs of type $3_n$ for prime $\frac{n}{4}$ is a non-decreasing function.
% going through all integers from $1$ to $19$.
\end{remark}

\begin{figure}
\begin{minipage}{7cm}
\centering
\epsfxsize=30mm
\epsfbox{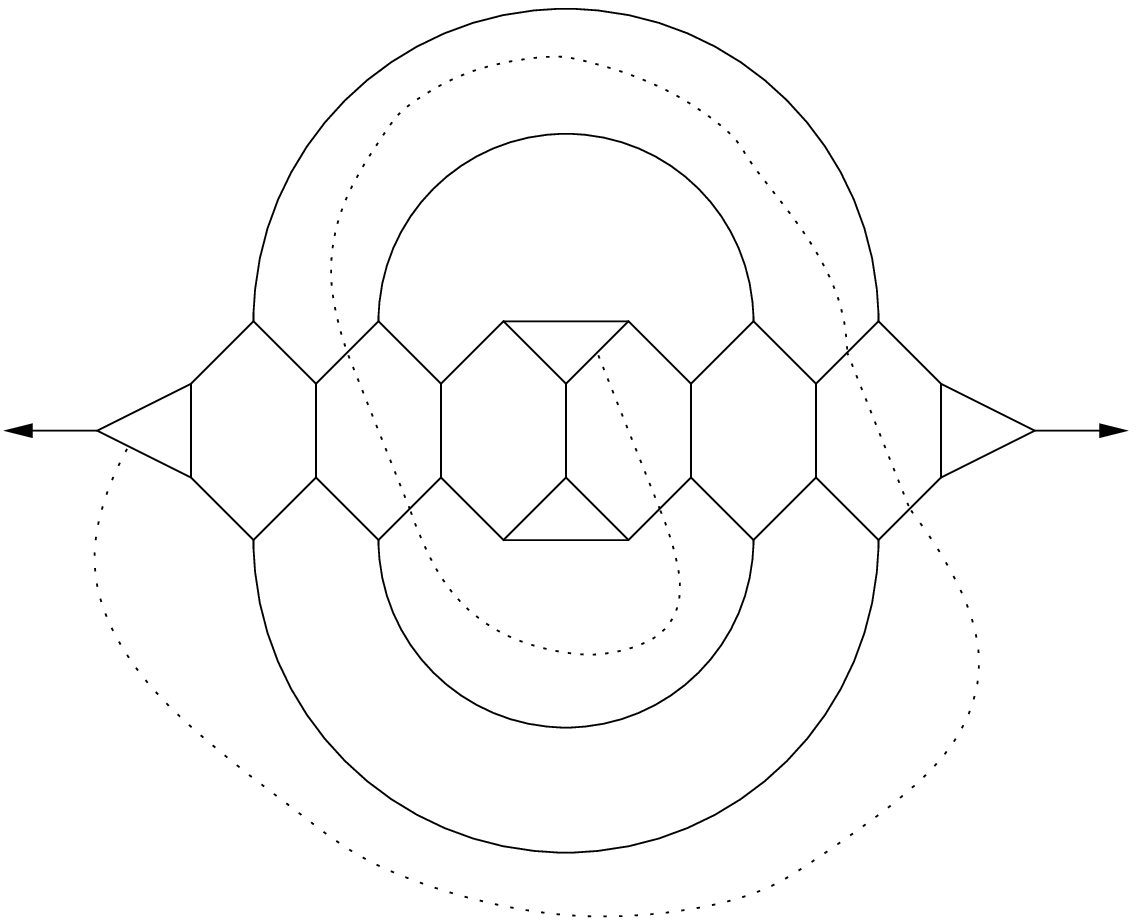}\par
a)
\end{minipage}
\begin{minipage}{7cm}
\centering
\epsfxsize=60mm
\epsfbox{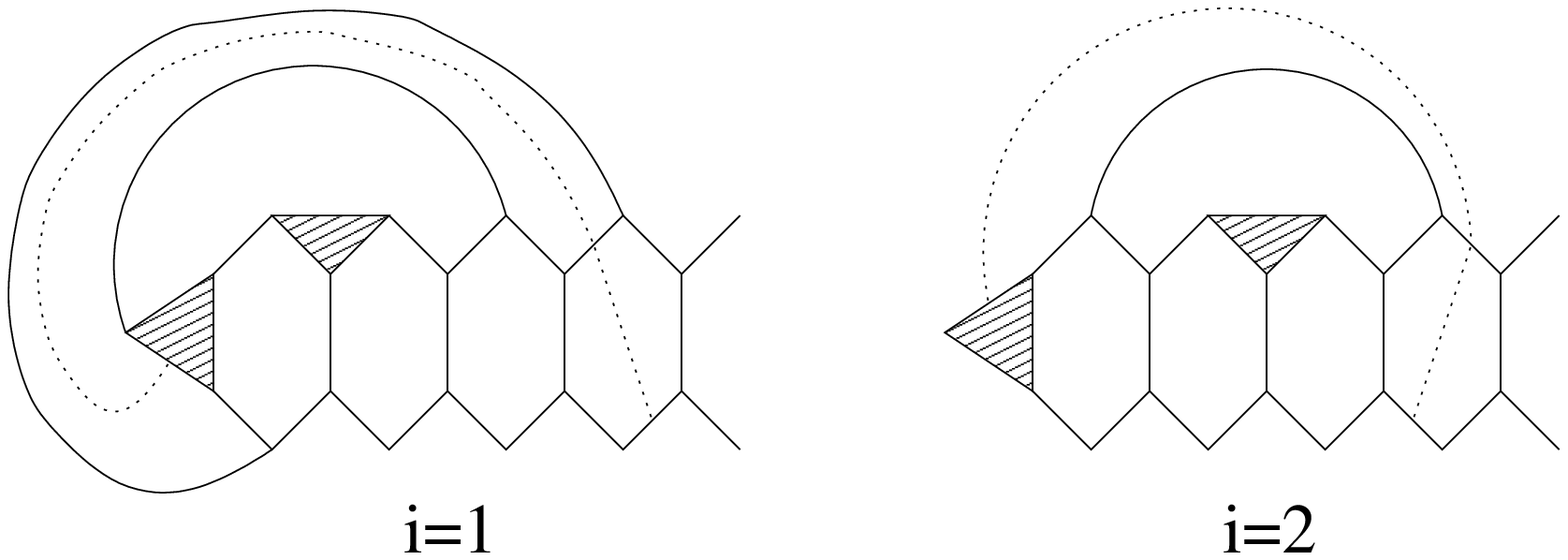}\par
b)
\end{minipage}
\caption{The pseudo-road construction.}
\label{fig:PseudoConstruction}
\end{figure}

\begin{theorem}
There exists a tight graph $G$ of type $3_n$ if and only if $\frac{n}{4}$ is odd.
%The class $3_n$ contains a tight graph if and only if $\frac{n}{4}$ is odd.
\end{theorem}
\proof The unique (for every integer $\frac{n}{4}\ge 4$) graph, defined in Theorem \ref{Possible-forms-for-3nS} (iii), is tight if $\frac{n}{4}$ is odd; we represent this graph in Figure \ref{fig:PseudoConstruction}.a) for $n=28$.

%\begin{center}
%\epsfxsize=60mm
%\epsfbox{T1T2T3T4.eps}
%\end{center}

Suppose now that $\frac{n}{4}$ is even and that $G$ is tight. We will use the necessary conditions of Theorem \ref{Theorem-interest-ZigZag-3n} to get a contradiction.

Consider the patch, formed by a pseudo-road $PR$ between two
triangles, say, $T_1$ and $T_2$.
Since $G$ is assumed to be tight, there are $\frac{n}{4}-1$ hexagons
in $PR$. Moreover, since there are no railroads, any of the
remaining two
triangles $T_3$ and $T_4$ will be adjacent to the hexagons of $PR$.
Let $i$ be the position of $T_3$ in $PR$.
The choice of position for $T_3$ and $T_4$ determines our graph of type $3_n$;
so, we need to show that every choice of $i$ leads to a non-tight graph
$3_n$. If we find a shorter pseudo-road, then we are
done.

Define a pseudo-road $PR'$ by starting with $T_1$ and taking the
upper edge if $i$ is even and the lower edge if $i$ is odd (see Figure \ref{fig:PseudoConstruction}.b)).
$PR'$ intersects with $PR$.
From the choice of $PR'$, the position 
$p$ of the first hexagon of intersection, of $PR'$ with $PR$,
satisfies $p\equiv 0\pmod 4$. Moreover, all positions of hexagons of
intersection satisfy this condition.
So, $PR'$ is shorter than $\frac{n}{4}-1$ and
the graph is {\em not} tight.\qed

The list of all symmetry groups of graphs $3_n$ is known (see
\cite{cremona}); they are (with their first appearance): $D_{2}(24_2)$, 
$D_{2h}(16_1)$, $D_{2d}(20_1)$, $T(28_2)$, $T_d(4_1)$ (see \cite{D1}
for the corresponding drawings).

%\begin{theorem}
%All symmetry groups of $3_n$ (with their first appearance) are $D_{2}(24_2)$, $D_{2h}(16_1)$, $D_{2d}(20_1)$, $T(28_2)$, $T_d(4_1)$.
%\end{theorem}
%\proof Let $G$ be a graph $3_n$ with triangles $T_i$, $1\leq i\leq 4$
%and symmetry group $Aut(G)$. The only symmetry, leaving 
%invariant all four
%triangles is the identity, since all triangles can not lie on a plane or
%on an axis of symmetry. So, the group $Aut(G)$ is isomorphic to a
%subgroup of the symmetry group of the Tetrahedron, which is $T_d$.
%
%Consider two triangles, say, $T_1$ and $T_2$, which are connected by a
%pseudo-road $PR$. The Gr\"{u}nbaum-Motzkin construction implies the
%existence of a pseudo-road, say, $PR'$, which connects the
%remaining two triangles $T_3$ and $T_4$.
%
%Define a rotation of order $2$, which exchange $T_1$ and $T_2$, by
%taking as its axis, the line connecting the middle of $PR$ (it can
%be an edge or an hexagon) and the middle of $PR'$. Another rotation
%of order $2$ is defined, which interchanges pseudo-road $PR$ and $PR'$.
%These two rotations generate a group of order four, which we can
%identify with the subgroup $D_2$ of $Aut(G)$. All subgroups of $T_d$,
%containing $D_2$, are $D_2$, $D_{2d}$, $D_{2h}$, $T$ and $T_d$, which
%are exactly all the known symmetry groups of $3_n$.\qed

For the point groups $D_{2}$, $D_{2d}$ and $D_{2h}$, we have the
following conjecture, which we checked for $n\leq 500$.

\begin{conjecture}
(i) $3_n(D_{2h})$ exists if and only if $\frac{n}{4}$ is even and $n\geq 16$; there are no tight $3_n(D_{2h})$.

(ii) $3_n(D_{2d})$ exists if and only if, either $\frac{n}{4}$ is odd and $n\geq 20$, or $\frac{n}{8}$ is even and $n\geq 24$; the graph defined in Theorem \ref{Possible-forms-for-3nS} (iii) is unique $3_n(D_{2d})$ tight graph if $\frac{n}{4}$ is odd and $n\geq 20$. 

(iii) $3_n(D_2)$ exists if and only if $n\geq 24$ and $n\not= 28,
32$; tight $3_n(D_2)$ exists if and only if $\frac{n}{4}$ is odd,
$n\geq 44$ and $n\not= 60, 84$; there are $i$ tight $3_n(D_2)$
starting with $\frac{n}{4}=6i+5$ for $1\leq i \leq 9$, except the case
$i=5$ starting with $\frac{n}{4}=37$.

\end{conjecture}

\section{Polyhedra $4_n$}
Since any graph of type $4_n$ is bipartite, we can apply
Theorem \ref{Shtogrin-Shank}: all self-intersections of every zigzag
are of type I. Moreover, there exists an orientation of zigzags,
with respect to which each edge has type I. We will always work
with such orientation. Any self-intersection, double or triple,
of railroad is of type I (see Figure \ref{fig:TypeSelfIntersection}).

We present in Table \ref{tab:listOf4nGraphs} the numbers $N_4(n)$ of
$4_n$ and the numbers $N^t_4(n)$ of tight $4_n$ for $n\leq 260$.

\begin{table}
\scriptsize
\begin{equation*}
\begin{array}{||c|c|c||c|c|c||c|c|c||c|c|c||c|c|c||c|c|c||}
\hline
\hline
n&N_4&N_4^{t}&n&N_4&N_4^{t}&n&N_4&N_4^{t}&n&N_4&N_4^{t}&n&N_4&N_4^{t}&n&N_4&N_4^{t}\\\hline
8&1&1&52&13&10&94&46&27&136&165&104&178&291&176&220&646&403\\
12&1&1&54&10&6&96&59&45&138&110&89&180&298&215&222&426&276\\
14&1&0&56&23&16&98&65&38&140&220&134&182&356&190&224&776&469\\
16&1&1&58&12&8&100&70&50&142&150&96&184&388&246&226&584&361\\
18&1&1&60&19&15&102&48&34&144&164&116&186&259&174&228&567&413\\
20&3&2&62&21&11&104&99&62&146&189&109&188&479&285&230&684&357\\
22&1&1&64&22&17&106&65&43&148&207&145&190&352&203&232&754&481\\
24&3&3&66&16&12&108&79&60&150&142&106&192&352&264&234&489&331\\
26&3&1&68&36&24&110&89&53&152&265&167&194&418&230&236&894&517\\
28&3&3&70&21&12&112&97&63&154&190&109&196&463&310&238&681&421\\
30&2&2&72&29&23&114&68&47&156&202&146&198&303&207&240&663&471\\
32&8&5&74&31&16&116&133&80&158&237&121&200&559&338&242&781&463\\
34&3&2&76&34&24&118&90&58&160&262&162&202&419&274&244&867&533\\
36&7&5&78&24&14&120&99&75&162&175&123&204&424&313&246&566&381\\
38&7&4&80&53&32&122&115&67&164&330&205&206&498&260&248&1016&594\\
40&7&5&82&32&21&124&127&81&166&239&140&208&553&355&250&790&466\\
42&5&4&84&42&33&126&86&58&168&249&193&210&365&251&252&751&541\\
44&14&7&86&47&25&128&171&103&170&288&154&212&663&400&254&902&487\\
46&6&5&88&50&36&130&118&79&172&319&221&214&500&309&256&1006&629\\
48&12&9&90&35&27&132&133&100&174&209&135&216&483&358&258&661&459\\
50&12&8&92&75&48&134&152&82&176&397&250&218&580&332&260&1173&675\\
\hline
\hline
\end{array}
\end{equation*}
\caption{Numbers $N_4(n)$ of graphs $4_n$ and $N_4^t(n)$ of tight $4_n$ for $n\leq 260$.}
\label{tab:listOf4nGraphs}
\end{table}

\begin{conjecture}

(i) $z$-knotted $4_n$ exists if and only if $n\geq 30$, $n\equiv 2 \pmod 4$.

(ii) Tight $4_n$ exists if and only if $n\geq 8$, $n\not= 10, 14$.

(iii) Every tight $4_n$ has at most eight zigzags.

\end{conjecture}

\begin{theorem}
All symmetry groups of $4_n$ (with their first appearance) are 
$C_1(40_4)$, $C_s(34_2)$, $C_2(26_1)$, $C_{i}(140_{12})$, $C_{2v}(22_1)$,
$C_{2h}(44_7)$, $D_2(24_2)$, $D_3(20_1)$, $D_{2d}(16_1)$, $D_{2h}(20_2)$,
$D_{3d}(20_3)$, $D_{3h}(14_1)$, $D_6(84_{42})$, $D_{6h}(12_1)$, 
$O(56_{17})$, $O_h(8_1)$.
\end{theorem}
\proof Let $G$ be a graph of type $4_n$ and let $r$ be a rotation
of $G$ of order $k$. By definition of a graph $4_n$, one has
$k=2$, $3$, $4$ or $6$.

If $k=6$, then the axis of $r$ goes through two hexagons,
say, $H_1$ and $H_2$. Consider around $H_1$ the ring of hexagons adjacent to it; then, after adding $p$ such concentric rings of hexagons, one will encounter a square and so, by $6$-fold symmetry, all squares of $G$.
One can then complete, in a unique way, the graph by adding
the same number of rings of hexagons on the other side of the sphere. This construction implies the existence of rotations of order two passing through opposite squares and middles of consecutive squares. This implies symmetry $D_6$ or $D_{6h}$.

If $k=4$, then the axis of $r$ goes through two squares,
say, $sq_1$ and $sq_2$. After adding $p$ rings of hexagons 
around $sq_1$, one finds a square and so, by symmetry, four
squares, say, $sq_3$, \dots, $sq_6$. One can complete the graph in 
a unique way; from the construction it is clear that there
is a $4$-fold axis through $sq_3$. So, the group is $O$ or $O_h$.

%If $k=3$, then the axis of $r$ goes, either through two vertices, two hexagons, or one vertex and one hexagon, which we denote by $C_1$, $C_2$.
%After adding $p$ rings of hexagons around $C_1$, one finds three squares; then, adding $q$ rings of hexagons one finds the three other squares.
%We complete the graph in, with hexagons, in an unique way; then, we remark that $C_1$ and $C_2$ are both hexagons or vertices and that the structure of $G$ around $C_1$ is identical to the one around $C_2$, so there are three rotations of order $2$ exchanging them and the symmetry group is $D_3$, $D_{3h}$, or $D_{3d}$.

If $k=3$, then the axis of $r$ goes, either through one vertex, or
one hexagon. After adding $p$ rings of hexagons around this center, 
one finds three squares; then, adding $q$ rings of hexagons, one finds 
the three other squares. The patch, formed by those six squares and
the $q$ rings of hexagons, has an additional symmetry, formed by 
exchanging the two triples of squares. So, $G$ has this symmetry too
and the symmetry group is $D_3$, $D_{3h}$, or $D_{3d}$.

Assume that $G$ has symmetry $S_4$. The six squares are partitioned in one orbit of two squares, say, $sq_1$ and $sq_2$, and another of four squares.
The $2$-fold axis goes through $sq_1$ and $sq_2$; moreover, if one consider the four pseudo-roads $PR_i$ from $sq_1$, then all of them stop at $sq_2$.
Consider the patch $P$ formed by $PR_1$ and $PR_2$. This patch has two acute angles; so, by Theorem \ref{Local-Euler-Formula}, it contains exactly one square $sq$.
We will show that $sq$ is in the plane, defined by $sq_1$ and $sq_2$, which will prove that there is a plane of symmetry for $G$. If $sq$ is next to $PR_1$ and $PR_2$, then their length is two and the result holds; otherwise, by removing some hexagons, we find a smaller patch with the same property and we conclude by induction.

The above reasonings yield the following list of possibilities: 
$C_1$, $C_s$, $C_2$, $C_{i}$, $C_{2v}$, $C_{2h}$, $D_2$, 
$D_3$, $D_{2d}$, $D_{2h}$, $D_{3d}$, $D_{3h}$, $D_6$, $D_{6h}$, $O$, $O_h$.
For every such group, there is a graph of type $4_n$ with $n\leq 140$
having this symmetry; their pictures are available from \cite{D1}. \qed

A similar study of point groups for $4$-valent plane graphs was done in \cite{DDS}.

\begin{figure}
\centering
\begin{minipage}[b]{6cm}%
\centering
\epsfig{figure=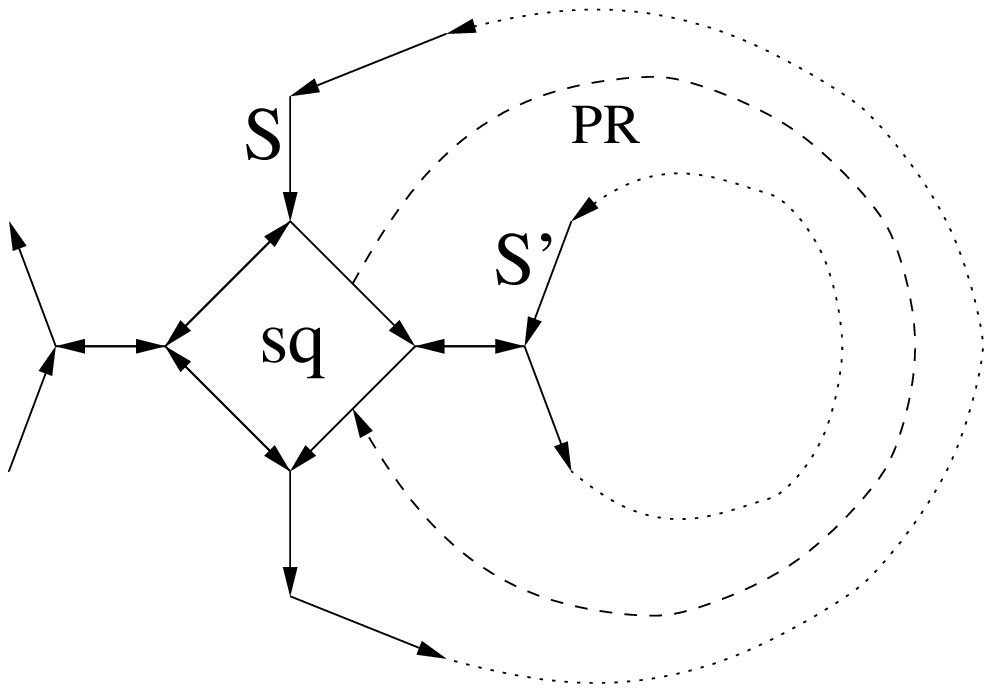,width=4cm}\par 
a) One lonely side $S$
\end{minipage}
\begin{minipage}[b]{6cm}%
\centering
\epsfig{figure=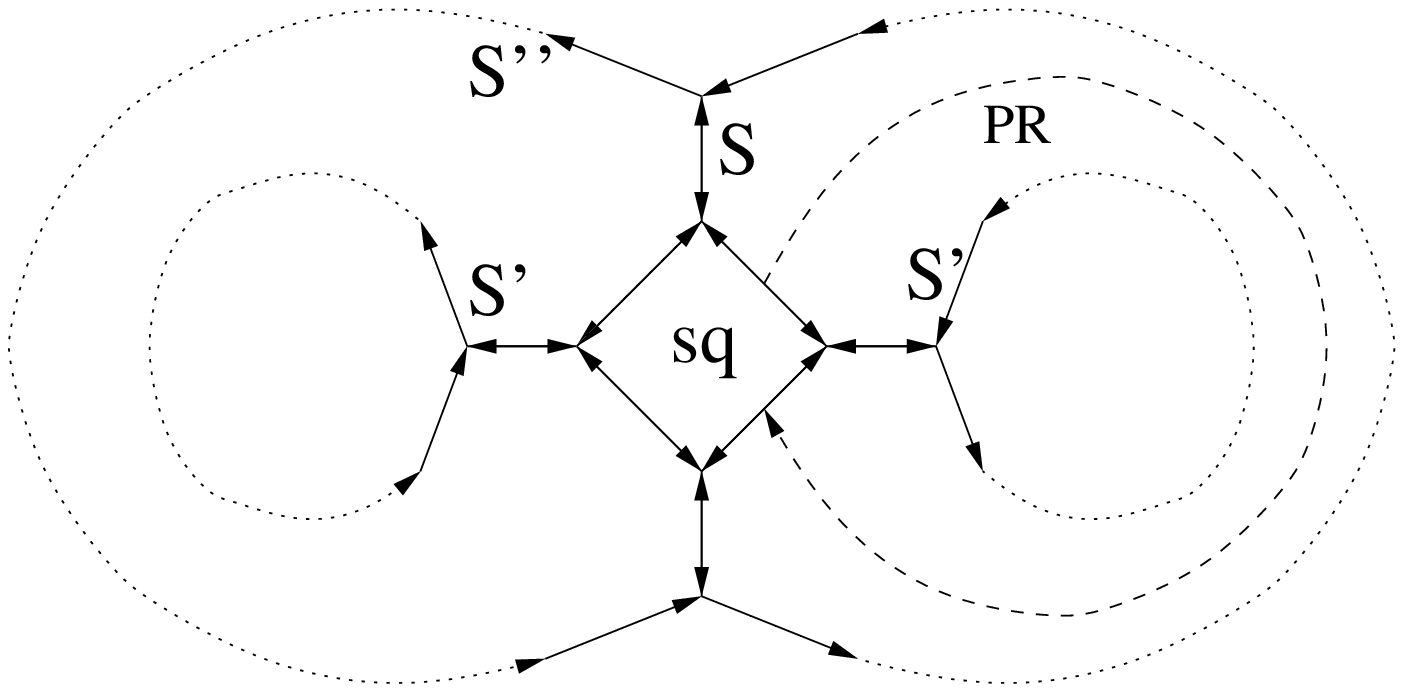,width=5.5cm}\par 
b) Two lonely sides $S$ and $S''$
\end{minipage}
\caption{Two cases of lonely sides.}
\label{fig:CaseLonely4n}
\end{figure}

\begin{theorem}\label{theoremUPPERBOUIND9}
Every tight $4_n$ has at most $9$ zigzags.
\end{theorem}
\proof Let $G$ be a tight $4_n$ with zigzags $Z_1$, \dots, $Z_l$.
%By Proposition \ref{crude-upper-bound}, $l\leq 12$.
We obtain $2l$ sides, since every zigzag has two sides.
A side $S$ is called
{\em lonely} if it is incident to exactly one square $sq$. We can define
a pseudo-road $PR$, parallel to $S$, which will begin and finish on $sq$.
On the other side of the pseudo-road, there will be a side $S'$ of a zigzag
$Z'$, which will be incident two times to $sq$. Moreover, if $S'$ is
incident exactly two times to $sq$, then it defines another lonely
side $S''$ (see Figure \ref{fig:CaseLonely4n}.b)).

Call $n_{1a}$ the number of lonely sides in first case and $n_{1b}$
the number of lonely sides in second case. Also call $n_2$ the number of 
sides incident to exactly two squares (identical or not) and $n_3$ the
number of sides incident to at least three squares (identical or not).
Obviously, $l=\frac{1}{2}(n_{1a}+n_{1b}+n_2+n_3)$.

In case a) of lonely side, $S'$ is incident to at least three
squares; so, $n_{1a}\leq n_3$, while in case b) $S'$ is incident exactly
two times to $sq$; so, $n_{1b}\leq n_2$. Every square $sq$ can be incident
to $0$, $1$ or $2$ lonely sides; so, one has the inequality 
$n_{1a}+\frac{n_{1b}}{2}\leq 6$. By an enumeration of incidences, one
gets $n_{1a}+n_{1b}+2n_2+3n_3\leq 4\times 6=24$. The $4$-dimensional
polyhedron defined by above inequalities and non-negativity inequalities
$n_{1a}, \dots, n_3\geq 0$ is denoted by ${\cal P}$. We can maximize the
linear function $\frac{1}{2}(n_{1a}+n_{1b}+n_2+n_3)$ over $P$ using the
cdd program (see \cite{cdd}) and obtain the upper bound $9$, which is 
attained for the unique $4$-uple $(n_{1a}, n_{1b}, n_2, n_3)=(0,12,6,0)$. \qed

The following theorem is in sharp contrast with Theorem
\ref{Theorem-interest-ZigZag-3n}(iv) for graphs $3_n$ and Theorem
\ref{TheoremBigDutourOne} for graphs $5_n$.

\begin{theorem}\label{TheoremIntersectionTwoSimpleZZ}
The intersection of every two simple zigzags of a graph of type $4_n$, if non-empty, has one of the following forms (and so, its size is $2$, $4$ or $6$).
\begin{center}
\epsfxsize=110mm
\epsfbox{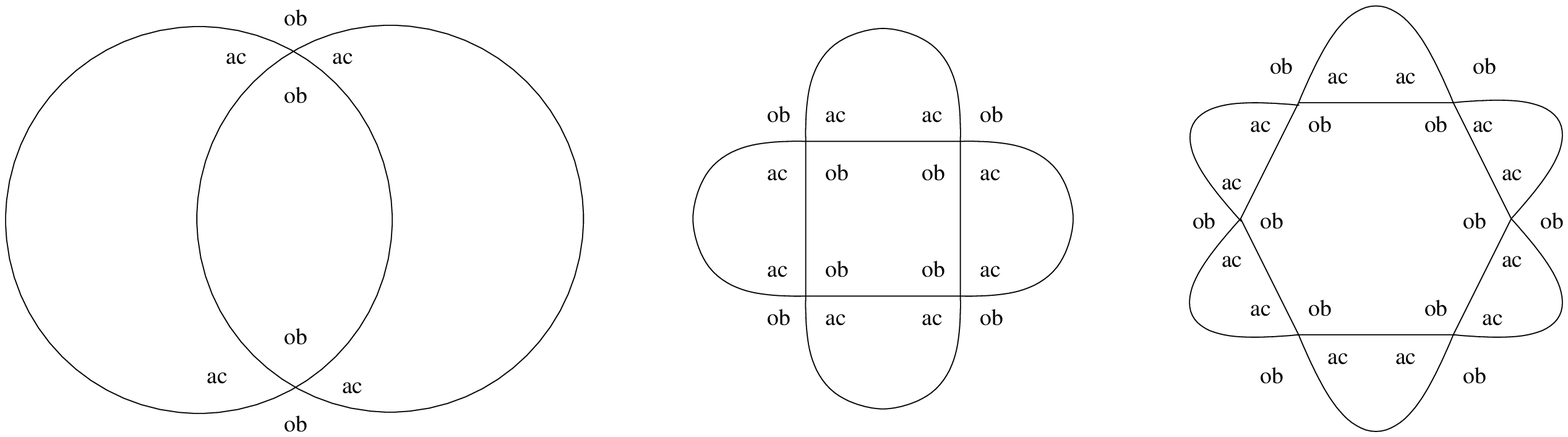}
\end{center}
\end{theorem}
\proof Define the graph $H$ as the graph, whose vertices are edges of
intersection between simple zigzags $Z$ and $Z'$, with two vertices
being adjacent if they are linked by a path belonging to one of $Z$,
$Z'$. The graph $H$ is a plane $4$-valent graph and $Z$, $Z'$ 
define two
central circuits in $H$. Since $Z$ and $Z'$ are simple,
the faces of $H$ are $t$-gons with even $t$.

Applying Theorem \ref{Local-Euler-Formula} to a $t$-gonal face $F$ of $H$, we obtain that the number $p'_4$ of $4$-gons in $F$ satisfies $6-t_{ob}-2t_{ac}=2p'_4$ .
So, the numbers $t_{ob}$ and $t_{ac}$ are even, since
$t=t_{ob}+t_{ac}$. Also, $6-t_{ob}-2t_{ac}\geq 0$. So, $t\leq 6$.

We obtain the following five possibilities for the faces of $H$: 
$2$-gons with two acute angles, $2$-gons with two obtuse angles, $4$-gons with four obtuse angles, $4$-gons with two acute and two obtuse angles, $6$-gons with six obtuse angles.

\begin{figure}
\begin{center}
\epsfxsize=110mm
\epsfbox{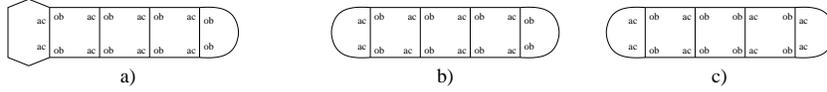}
\end{center}
\caption{Three cases for sequence of $4$-gons.}
\label{fig:CaseSequence}
\end{figure}

Take an edge $e$ of a $6$-gon in $H$ and consider the sequence
(possibly, empty) of adjacent $4$-gons of $H$ emanating from this
edge. This sequence will stop at a $2$-gon or a $6$-gon; the
case-by-case analysis of angles yields that this sequence has to stop
at a $2$-gon (see Figure \ref{fig:CaseSequence}.a)).

Take an edge of a $2$-gon in $H$ and consider the same construction. 
If the angles are both obtuse, then the construction is identical and 
the sequence will terminate at a $2$-gon or a $6$-gon. If the angles are 
both acute, then cases b), c) of Figure \ref{fig:CaseSequence} are possible.

In the first case, all $4$-gons contain two obtuse angles and two acute
angles; so, the pseudo-road finishes with an edge of two obtuse
angles. In the second case, there is a $4$-gon, whose angles are all
obtuse; this $4$-gon is unique in the sequence and its position is
arbitrary. Every pair of opposite edges of a $4$-gon belongs to a 
sequence of $4$-gons considered above. So, all angles of a $4$-gon 
are the same, i.e., obtuse. This fact restricts 
the possibilities of intersections to the three cases of the theorem. \qed

\begin{corollary}
The only tight $4_n$, having only simple zigzags, are the Cube and the Truncated Octahedron.
\end{corollary}
\proof By Theorem \ref{TheoremIntersectionTwoSimpleZZ}, every two
simple zigzags intersect in at most six edges. Since, by Theorem
\ref{theoremUPPERBOUIND9}, there are at most $9$ zigzags, we obtain
the upper bound $(9-1)6=48$ on the length of every zigzag. This yields the
upper bounds $\frac{9}{2}48=216$ on the number of edges of $G$
and $\frac{2}{3}216=144$ on the number of its vertices. But, our
exhaustive computation in this range
of values gave the Cube and the Truncated Octahedron as the only
solutions. \qed

\section{Polyhedra $5_n$}
Zigzag structure of fullerenes was studied in \cite{DDF}, for which
present paper is a follow-up. Here, we add only several observations.

The smallest $z$-unbalanced fullerenes $5_n$ with only simple zigzags are $5_{108}(D_{2d})$ and $5_{144}(D_3)$  with $z$-vectors $24^{8}$, $26^{4}$, $28$ and $28^{12}$, $32^{3}$, respectively. Any $z$-uniform or tight $5_n$, having only simple zigzags, with $n\leq 200$ is $z$-balanced.

Table \ref{tab:t3} gives all tight $5_n$, $n\leq 200$, with
simple zigzags; we conjecture that this list is complete. Apropos,
amongst the $9$ fullerenes of Table \ref{tab:t3} only $5_{60}(I_h)$,
$5_{88}(T)$, $5_{140}(I)$ have $12$ isolated pentagons, and only
$5_{76}(D_{2d})$ has $6$ isolated pairs of adjacent pentagons.

\begin{table}
\scriptsize
\begin{center}
\begin{tabular}{||c|c|l|l|c||}
\hline
\hline
n       &Group          &$z$-vector       &Orb. sizes     &Int. vector\\
\hline \hline
$20$    &$I_h$          &$10^6$         &6              &$2^5$ \\
$28$    &$T_d$          &$12^7$         &3,4            &$2^6$\\
$48$    &$D_3$          &$16^9$         &3,3,3          &$2^8$\\
$60$    &$I_h$          &$18^{10}$              &10             &$2^9$\\
$60$    &$D_3$          &$18^{10}$              &1,3,6          &$2^9$\\
$76$    &$D_{2d}$       &$22^4,20^7$            &1,2,4,4        &$4,2^9$ and $2^{10}$\\
$88$    &$T$            &$22^{12}$              &12             &$2^{11}$\\
$92$    &$T_d$          &$24^4,20^9$            &3,4,6          &$2^{12}$ and $2^{10},0^2$\\
$140$   &$I$            &$28^{15}$              &15             &$2^{14}$\\
\hline
\hline
\end{tabular}
\end{center}
\caption{All tight fullerenes $5_n$ with $n\leq 200$, having only simple zigzags.}
\label{tab:t3}
\end{table}

It was conjectured in \cite{DDF} that for any even $n\geq 20$,
$n\not=22$, there exists a tight $5_n$. Moreover, we could not find a
tight graph $5_n$ with more than $15$ zigzags.

The list of all symmetry groups of fullerenes $5_n$ is known (see
\cite{FM}); they are (with their first appearance): $C_1(36_3)$,
$C_2(32_1)$, $C_i(56_{314})$, $C_s(34_2)$, $C_3(40_{30})$, $D_2(28_1)$, 
$S_4(44_{82})$, $C_{2v}(30_2)$, $C_{2h}(48_{80})$, $D_3(32_6)$, 
$S_6(68_{6263})$, $C_{3v}(34_6)$, $C_{3h}(62_{2334})$, $D_{2h}(40_{33})$,
$D_{2d}(36_6)$, $D_5(60_{1794})$, $D_6(72_{11144})$, $D_{3h}(26_1)$, 
$D_{3d}(32_3)$, $T(44_{73})$, $D_{5h}(30_1)$, $D_{5d}(40_1)$, 
$D_{6h}(36_{15})$, $D_{6d}(24_1)$, $T_d(28_2)$, $T_h(92_{126311})$, 
$I(140_{x})$, $I_h(20_1)$.

\begin{figure}
\centering
\begin{minipage}[b]{6.5cm}%
\leavevmode
\begin{center}
\epsfig{figure=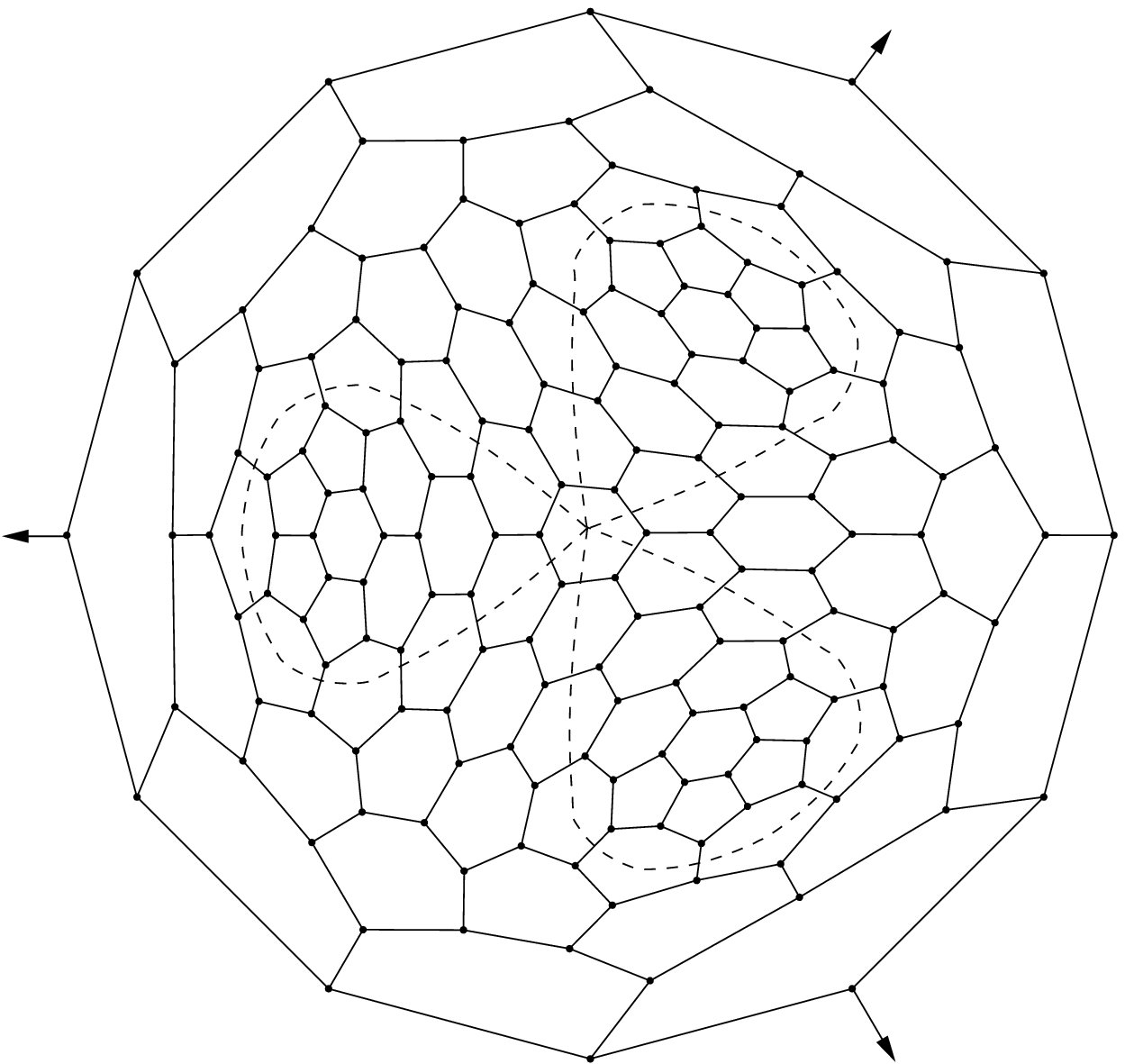,height=5cm}
\end{center}
\begin{center}
a) $5_{172}(C_{3v})$ with triple self-intersecting railroad
\end{center}
\end{minipage}
\begin{minipage}[b]{6.5cm}%
\leavevmode
\begin{center}
\epsfig{figure=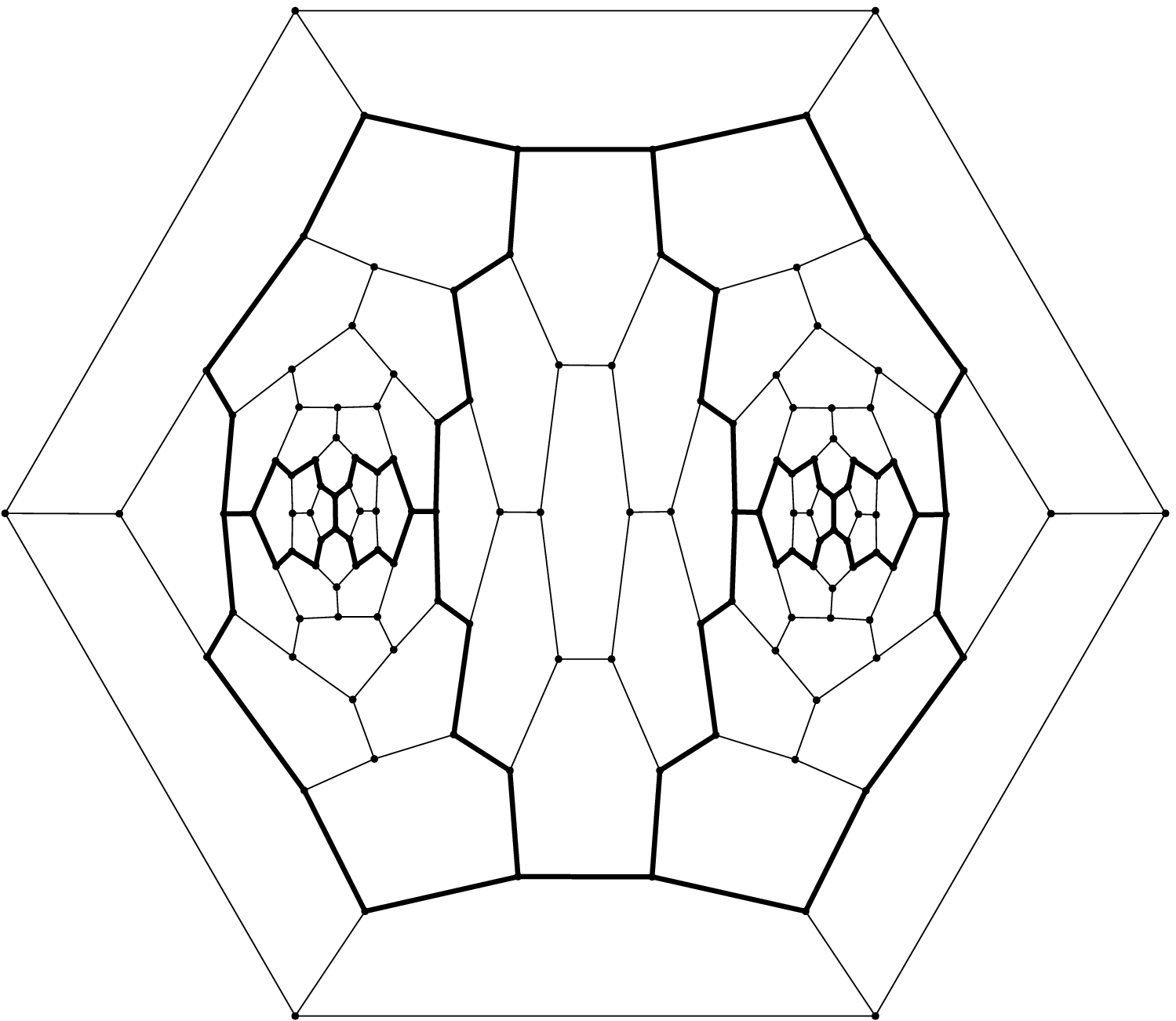,height=5cm}
\end{center}
\begin{center}
b) Two zigzags in $5_{136}(D_{2h})$ with $|Z\cap Z'|=8$
\end{center}
\end{minipage}
\caption{Two particular fullerenes.}
\label{fig:constructedByHand}
\end{figure}

There exist fullerenes admitting railroads with {\em triple}
self-intersection; see, for example, Figure \ref{fig:constructedByHand}.a).
More generally, we believe that any curve, which represents 
a railroad in a graph of type $4_n$, also appears as a curve
representing a railroad in some $5_{n'}$.

\begin{theorem}\label{TheoremBigDutourOne}
For any even number $h\geq 2$, there exists a fullerene $5_n$, 
with $n=18h-8$, having two simple zigzags intersecting 
in exactly $h$ edges. It has symmetry $T_d$ if $h=2$ and 
$D_{2h}$, $D_{2d}$ if $\frac{h}{2}$ is even, odd respectively.
\end{theorem}
\proof For any even $h\geq 2$ there exist a unique $h$-vertex
$4$-valent plane graph $H$, whose faces are four $2$-gons
(in two pairs of adjacent ones) and $\frac{h}{2}-2$ $4$-gons only, 
and having two simple central circuits (see \cite{DS}). This graph
has symmetry $D_{4h}$ for $h=2, 4$ and, for larger values, symmetry
$D_{2h}$, $D_{2d}$ if $\frac{h}{2}$ is even, odd respectively.

We will identify each of the two central circuits of $H$ with simple
zigzags, $Z_1$ and $Z_2$, and each vertex with an edge of intersection
between them. Every face of $H$ can be seen as a patch in the graph
$5_n$ which we will construct, and so, the local Euler formula
(\ref{Local-Euler-Formula}) can be applied. Fix a face $F$ of $H$; 
one can assign to every angle of
$F$ an angle (obtuse or acute), so that every $2$-gon has one acute
and one obtuse angle, while every $4$-gon has two obtuse and two
acute angles. See below the graph for the first values of $h$.

\begin{center}
\epsfxsize=90mm
\epsfbox{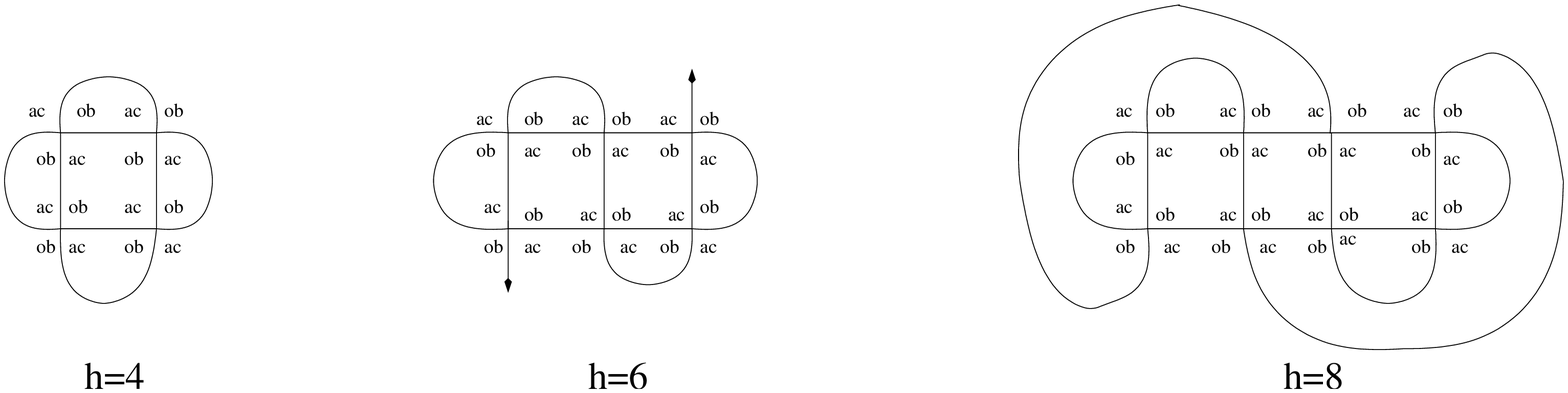}
\end{center}

So, $2$-gonal patches will contain three pentagons, while
$4$-gonal patches will contain only hexagons. We replace $2$-
and $4$-gonal faces of $H$ by patches depicted below.

\begin{center}
\epsfxsize=90mm
\epsfbox{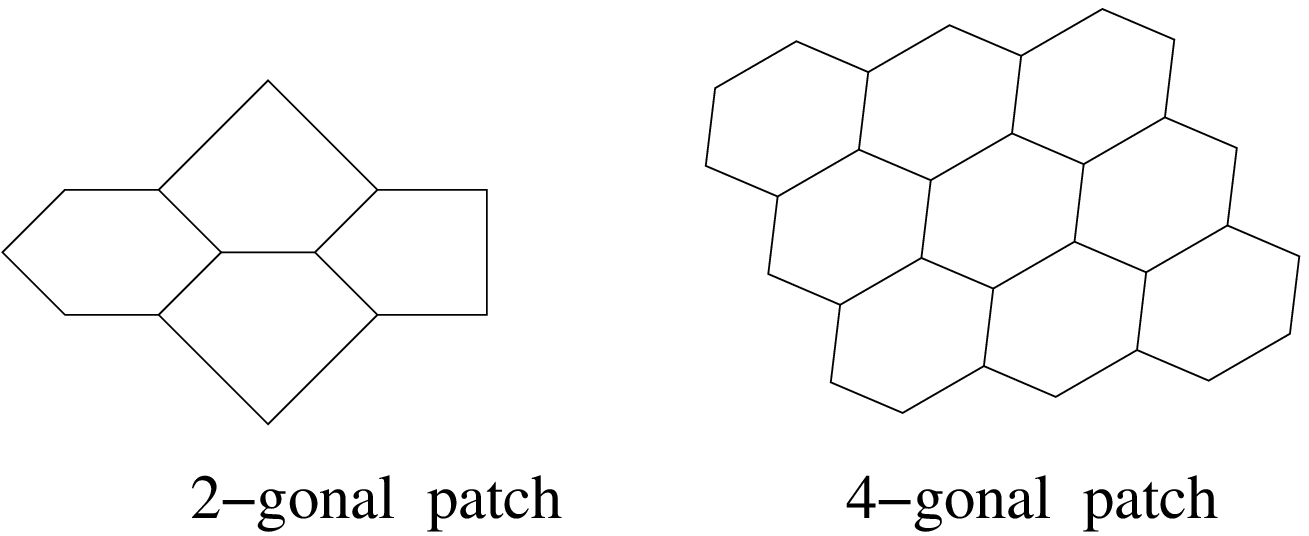}
\end{center}

The obtained graph has $9(h-2)+4$ hexagonal faces, and so,
$n=18h-8$ vertices. The symmetry group is the same as the 
one of $H$, except for the first values $h=2,4$. See
Figure \ref{fig:constructedByHand}.b) for the corresponding graph with
$h=8$.\qed

%\begin{center}
%\epsfxsize=30mm
%\epsfbox{ThreePentagons.eps}
%\end{center}

\end{document}